\documentclass[12pt,reqno]{amsart}
\usepackage{amscd} 
\usepackage{fancybox}

\usepackage[hyphens]{xurl}
\usepackage{mathrsfs}
\usepackage{amsthm}
\usepackage{amsmath}
\usepackage{amssymb}
\usepackage{amsfonts}
\usepackage{pb-diagram} 
\usepackage[all]{xy}
\usepackage[dvipdfmx]{graphicx}
\usepackage{latexsym}
\usepackage{comment}
\usepackage{ascmac}
\usepackage{caption}
\usepackage{mathtools}
\usepackage
{lineno}
\usepackage{color}
\def\qed{\hfill $\Box$}

\numberwithin{equation}{section}

\newcommand{\diam}{{\rm diam}}
\newcommand{\dimH}{\dim_{\rm H}}

\makeatletter

    \newcommand\contFrac{\@ifstar{\@contFracStar}{\@contFracNoStar}}

   \def\singleContFrac#1#2{%
        \begin{array}{@{}c@{}}%
            \multicolumn{1}{c|}{#1}%
            \\%
            \hline%
           \multicolumn{1}{|c}{#2}%
        \end{array}%
   }

    \def\@contFracNoStar#1{%
        \mathchoice{
            \@contFracNoStarDisplay@#1//\@nil%
        }{
            \@contFracNoStarInline@#1//\@nil%
        }{
            \@contFracNoStarInline@#1//\@nil%
        }{
            \@contFracNoStarInline@#1//\@nil%
        }%
    }

    \def\@contFracNoStarDisplay@#1//#2\@nil{%
        \@ifmtarg{#2}{%
            #1%
        }{%
            #1+\cfrac{1}{\@contFracNoStarDisplay@#2\@nil}%
        }%
    }

        \def\@contFracNoStarInline@#1//#2\@nil{%
            \@ifmtarg{#2}{%
                #1%
            }{%
                #1 \@@contFracNoStarInline@@#2\@nil%
            }%
        }
        \def\@@contFracNoStarInline@@#1//#2\@nil{%
            \@ifmtarg{#2}{%
                + \singleContFrac{1}{#1}%
            }{%
                + \singleContFrac{1}{#1} \@@contFracNoStarInline@@#2\@nil%
            }%
        }

    \def\@contFracStar#1{%
        \mathchoice{
            \@contFracStarDisplay@#1////\@nil%
        }{
            \@contFracStarInline@#1//\@nil%
        }{
            \@contFracStarInline@#1//\@nil%
        }{
            \@contFracStarInline@#1//\@nil%
        }%
    }

    \def\@contFracStarDisplay@#1//#2//#3\@nil{%
        \@ifmtarg{#2}{%
            #1%
        }{%
            #1 + \cfrac{#2}{\@contFracStarDisplay@#3\@nil}%
        }%
    }

        \def\@contFracStarInline@#1//#2\@nil{%
            \@ifmtarg{#2}{%
                #1%
            }{%
                #1 \@@contFracStarInline@@#2\@nil%
            }%
        }
        \def\@@contFracStarInline@@#1//#2//#3\@nil{%
            \@ifmtarg{#3}{%
                + \singleContFrac{#1}{#2}%
            }{%
                + \singleContFrac{#1}{#2} \@@contFracStarInline@@#3\@nil%
            }%
        }

\theoremstyle{plain}

\newtheorem{thm}{Theorem}[section]
\newtheorem{lem}[thm]{Lemma}
\newtheorem{cor}[thm]{Corollary}

\newtheorem{pro}[thm]{Proposition}
\theoremstyle{definition}

\newtheorem{rem}[thm]{Remark}

\newtheorem*{prf*}{Proof}

\newtheorem*{cla*}{Claim}
\newtheorem*{pf*}{}
\newtheorem*{lem*}{LemmaA}
\newtheorem*{lm*}{LemmaB}
\newtheorem*{stra*}{Strategy for the proof of main result A}
\makeatletter
\@namedef{subjclassname@2020}{%
  \textup{2020} Mathematics Subject Classification}
\makeatother
\title[Fractal transference principle for IFS]
{Fractal transference principle \\ for conformal iterated function systems}

\author{Yuto Nakajima and Hiroki Takahasi}

\usepackage{lineno}
\date{}

\address{Faculty of Science and Engineering, Doshisha University, Kyoto, 610-0394, JAPAN}
\email{yunakaji@mail.doshisha.ac.jp}
\address{Keio Institute of Pure and Applied Sciences (KiPAS), Department of Mathematics,
Keio University, Yokohama,
223-8522, JAPAN} 
\email{hiroki@math.keio.ac.jp}

\begin{document}

\begin{abstract} 

We establish a complex one-dimensional fractal transference principle connecting density combinatorics with the fractal geometry of complex continued fractions. Our framework applies to a broad class of conformal iterated function systems indexed by lattice points, including those associated with the Hurwitz continued fraction, and transfers positive-density combinatorial patterns to digit sets while preserving the maximal Hausdorff dimension.
As applications, we obtain fractal counterparts of the multidimensional and polynomial Szemerédi theorems, as well as their relativized versions yielding a counterpart of the multidimensional Green--Tao theorem.

\end{abstract}

\maketitle

\if0\begin{itemize}
\item Density Combinatorics: Studies the existence and properties of structured subsets (like arithmetic progressions) within subsets of natural numbers that have positive upper density (meaning they are sufficiently "thick"). 
\item Fractal Dimension Theory: Focuses on the "fractal dimension" (a measure of complexity) of sets of numbers, particularly those related to their continued fraction representations. 
\item The Bridge: The "fractal transference principle" builds a connection between these two areas, providing a way to translate combinatorial results from the realm of sequences of integers to the geometric context of real numbers and their continued fraction properties. \end{itemize}
\fi


\section{Introduction}
Density combinatorics
 aims to study the existence of patterns within subsets of integers that have positive density. The dimension theory of continued fractions studies the fractal dimensions of sets defined by conditions on continued fraction expansions. The {\it fractal transference principle}, coined in \cite{NT}, builds a connection between these two areas, providing a way to transfer density combinatorics theorems to the geometric context of real and complex numbers. 

From \cite{NT}, let us recall how the principle works for the regular continued fraction
\[
x=
\cfrac{1}{
a_1(x)
+\cfrac{1}{
a_2(x)
+\cfrac{1}{
a_3(x)
+\ddots
}
}
},
\]
where $x\in(0,1)$ is an irrational number and $a_n(x)\in\mathbb N$ for  $n\geq1$. Consider the set
\[E=\{x\in (0,1)\setminus\mathbb Q\colon a_m(x)\neq a_{n}(x)\ \text{ for all } m,n\in\mathbb N\text{ with }m\neq n\},\]
which is of Hausdorff dimension $1/2$ \cite{R85}.
For each $x$ in this set, $\{a_n(x)\colon n\in\mathbb N\}$ is naturally regarded as an infinite subset of $\mathbb{N}$. 
Let ($*$) be a certain property that holds for any subset of $\mathbb N$ with positive density, such as the existence of arbitrarily long arithmetic progressions \cite{S75} or polynomial progressions \cite{BL}.
 The principle asserts that for any $S\subset\mathbb N$ with positive density, there is a set $E_S\subset E$
 with Hausdorff dimension $1/2$
such that for any $x\in E_S$ the set  $\{a_n(x)\colon n\in\mathbb N\}\cap S$ inherits property ($*$). This yields ``fractal counterparts'' of  density combinatorics theorems such as Szemer\'edi's theorem~\cite{S75} and the Bergelson-Leibman theorem \cite{BL}, see \cite[Theorem~1.3]{NT} for more details.

Related developments have pursued the connection between additive patterns
and continued fractions in different directions. Tong and Wang~\cite{TW}
studied arithmetic progressions in the partial quotients of regular continued
fractions, and Zhang and Cao~\cite{ZC2} considered related positive-density
digit problems for infinite iterated function systems.
More recently, the first and second authors and Wang~\cite{NTW} studied
quantitative restrictions on arithmetic progressions in the partial quotients
of regular continued fractions, obtaining exact Hausdorff dimension formulas
depending on their prescribed lengths and positions.

Guo et al.  \cite{GOQZ} established a fractal
transference principle for real vectors by applying regular
continued fraction expansions to each coordinate.
The first author \cite{NakLaurent} established another extension
for continued fractions whose digits are polynomials over a finite
field.

The aim of this paper is to establish a complex one-dimensional
version of the fractal transference principle.
Our main target is complex continued fractions. To a given complex continued fraction one can often associate a natural iterated function system (IFS) consisting of infinitely many maps of the form
\[z\in\Delta\mapsto \frac{1}{z+\tau(i)},\ i\in I,\]
where $\Delta$ is a convex compact set in $\mathbb C$, $I$ is a countably infinite subset of $\mathbb Z^2$ and $\tau\colon I\to\mathbb C$ is an injection taking values in a subring of $\mathbb Q(\sqrt{-d})$, $d\in\mathbb N$.
In the examples considered below, points whose partial quotients tend
to infinity are, up to finitely many initial digits, represented by
limit sets of this type.
In this paper we provide a general, unified framework that is capable of dealing with a wide class of conformal IFSs and 
complex continued fractions.
In what follows, we first introduce the class of conformal iterated
function systems considered in this paper, then recall the
multidimensional density-combinatorial results that will be transferred
to their limit sets, and finally state our main results.

\subsection{Conformal IFS}
We use standard terminology from the theory of conformal iterated
function systems; see, for example, Mauldin and Urba\'nski~\cite{MauUrb}
for a comprehensive account.
Let $\Omega\subset \mathbb C$ be a region.
For a holomorphic function $\phi\colon\Omega\to\mathbb C$ and $x\in\Omega$ let $D\phi(x)$ denote the complex derivative of $\phi$ at $x$. 
Throughout this paper, let $\Delta\subset \mathbb C$ be a connected, convex compact set
such that
 the closure of its interior 
coincides with $\Delta$.   
Let 
$I$ be a countable subset of $\mathbb Z^2\setminus\{0\}$ 
and let $\{\phi_i\colon i\in I\}$ be a collection of maps from $\Delta$ to itself. For $n\geq2$ and $(i_1,\ldots,i_n)\in I^n$,
write \[\phi_{i_1\cdots i_n}=\phi_{i_1}\circ\cdots\circ\phi_{i_n}.\]
We say that $\{\phi_i\colon i\in I\}$ is a {\it conformal IFS} (iterated function system)  if the following three conditions hold:
\begin{itemize}

\item[(A1)] (open set condition) For
all $i,j\in I$ with $i\neq j$,
 \[\phi_{i}({\rm int}\Delta)\cap \phi_{j}({\rm int}\Delta)=\emptyset.\]

\item[(A2)] (conformality) There exists a connected open set $\widetilde \Delta\subset\mathbb C$ containing $\Delta$ 
such that each $\phi_i$ extends to a $C^{1}$ conformal diffeomorphism $\widetilde\phi_i\colon\widetilde \Delta\to\widetilde \phi_i(\widetilde \Delta)\subset\widetilde \Delta$.
With a slight abuse of notation, we often write $\phi_i$ for $\widetilde\phi_i$.

\item[(A3)](uniform contraction) There exist  $m\in\mathbb N$ and $\xi\in (0,1)$ such that for any $(i_1,\ldots, i_m)\in I^{m}$ and any $x\in \Delta$ we have 
\[0<|D\phi_{i_1\cdots i_m}(x)|\le \xi.\]
\end{itemize}
If, moreover,  $\#I=\infty$, then we say $\{\phi_i\colon i\in I\}$ is an {\it infinite conformal IFS} on $\Delta$.

Let $\Phi=\{\phi_i\colon i\in I\}$ be an infinite conformal IFS on $\Delta$. 
Fix $x_*\in\Delta$, and define an {\it address map} $\Pi\colon I^{\mathbb N}\to \Delta$ by 
\[\Pi((i_n)_{n=1}^{\infty})=\lim_{n\to \infty}\phi_{i_1\cdots i_n}(x_*).\]
The set $L(\Phi)=\Pi(I^{\mathbb N})$ is called the {\it limit set} of $\Phi$. By (A3), the address map is well-defined and does not depend on the choice of the point $x_*$.
Since we do not require ``$\phi_{i}(\Delta)\cap \phi_{j}(\Delta)=\emptyset$ for all $i,j\in I$ with $i\neq j$'', the address map may not be injective.
For each $x\in L(\Phi)$, any element of $\Pi^{-1}(x)$ is called an {\it address} of $x$.
Let $L'(\Phi)$ denote the set of $x\in L(\Phi)$ such that $\Pi^{-1}(x)$ is a singleton.
 For each $x\in L'(\Phi)$, 
let $(i_n(x))_{n=1}^\infty$ denote the address of $x$. 
Set \[
\mathcal O(\Phi):=
\bigcup_{\substack{i,j\in I\\i\neq j}}
\left(\phi_i(L(\Phi))\cap\phi_j(L(\Phi))\right).
\]
To ensure that $L(\Phi)\setminus L'(\Phi)$ is negligible, we 
introduce the following condition.
\begin{itemize}
\item[(A4)] 
$\mathcal O(\Phi)$ is of Hausdorff dimension $0$.
\end{itemize}
Condition (A4) implies that the set $L(\Phi)\setminus L^{\prime}(\Phi)$ 
has Hausdorff dimension zero (see Lemma~\ref{lem:negligible-overlaps}). 
\begin{rem}
The notion of a $2$-decaying conformal IFS was previously used in
\cite{NT2} under a stronger boundary condition.
The present condition {\rm (A4)} is a weaker formulation that
directly captures the negligibility of points with non-unique
addresses; see Lemma~\ref{lem:negligible-overlaps}.
\end{rem}

We introduce the maximum norm $\|x\|_\infty=\max\{|{\rm Re}(x)|,|{\rm Im}(x)|\}$ for $x\in \mathbb C$.
Note that $\|x\|_\infty\leq|x|\leq\sqrt{2}\|x\|_\infty$. 
We say $\Phi$ is {\it $2$-decaying} if (A4) holds, and there exists $C_0\geq1$ such that
for any $i\in I$ and any $x\in \Delta$ we have 
\begin{equation}\label{CF-der}\frac{C_0^{-1}}{\|i\|_\infty^{2 }}\le |D\phi_i(x)|\le \frac{C_0 }{\|i\|_\infty^{2} }.\end{equation}
More generally, one can consider $\nu$-decaying conformal IFSs $(\nu>1)$ replacing $\|i\|_\infty^2$ in \eqref{CF-der} by $\|i\|_\infty^\nu$.
The case $\nu=2$ is related to various complex continued fractions
and is therefore the most relevant one for the present paper.

\subsection{Multidimensional density combinatorics theorems}
Let $A\subset\mathbb Z^2$. A subset of $\mathbb Z^2$ of the form $v+nA$, $v\in\mathbb Z^2$, $n\in\mathbb N$ 
is called {\it a homothetic copy of $A$}.
As a multidimensional extension of Van der Waerden's theorem \cite{Wae27}, 
Gallai (see \cite[Theorem 8]{Gra90} and \cite[p.123]{R39}) proved that if
$r\in\mathbb N$ and
 $\mathbb Z^2$ is partitioned into $r$ subsets then 
 for every finite subset \(A\subset\mathbb Z^2\), one of them contains a homothetic copy of \(A\).

Gallai's theorem does not provide any information on which subset contains a homothetic copy of any finite subset of $\mathbb Z^2$.
A significant strengthening was obtained by Furstenberg and Katznelson \cite{FK78}, who proved that any subset
of $\mathbb Z^2$ with positive upper Banach density contains a homothetic copy of any finite subset of $\mathbb Z^2$.
For a set $S\subset \mathbb Z^2$ and
 a sequence of axis-parallel integer rectangles 
$\Pi_k=[p_k^{(1)},q_k^{(1)}]
\times [p_k^{(2)},q_k^{(2)}]$ in $\mathbb R^2$, $k\in\mathbb N$ with $q_k^{(j)}-p_k^{(j)}\to\infty$ as $k\to\infty$ for $j=1, 2$, define
\[d_{\rm B}(S,(\Pi_k)_{k=1}^\infty)=\limsup_{k\to\infty}\frac{\#(S\cap\Pi_k)}{\#(\mathbb Z^2\cap\Pi_k)}.\]
The {\it upper Banach density} of $S$, denoted by $d_{\rm B}(S)$, is given by
\[d_{\rm B}(S)=\sup d_{\rm B}(S,(\Pi_k)_{k=1}^\infty),\]
where the supremum is taken over all sequences of axis-parallel integer rectangles in $\mathbb R^2$.
The {\it upper density} of $S$, denoted by $\overline{d}(S)$,   
is given by \[\overline{d}(S)=\limsup_{N\to \infty}\frac{\#({S\cap[-N,N]^2})}{(2N+1)^2}.\]
Notice that $d_{\rm B}(S)\geq \overline{d}(S)$.


Bergelson and Leibman \cite{BL} proved a vast generalization of Furstenberg-Katznelson's theorem, known as polynomial Szemer\'edi's theorem. A precise statement of Bergelson-Leibman's theorem will be given later.

\subsection{Statements of the results}
For convenience, we identify $\mathbb C$ with $\mathbb R^2$. 
We restrict ourselves to the set 
\[E(\Phi)=\{x\in L'(\Phi)\colon i_{m}(x)\neq i_{n }(x)  \ \text{\rm for all } m,n\in\mathbb N\text{ with }  m\neq n\}.\]
Each $x\in E(\Phi)$  induces an injection 
$n\in \mathbb N\mapsto i_n(x)\in I$.
Therefore, 
$\{i_n(x)\colon n\in\mathbb N\}$ is naturally regarded as an infinite subset of $\mathbb{Z}^2$.
Let $\dim_{\rm H}$ denote the Hausdorff dimension on $\mathbb R^2$.

Before stating our main results, we introduce a growth condition on subsets of
$\mathbb{Z}^2$. For an infinite set $J\subset \mathbb{Z}^2$ and $N>0$, put
\[
Q_N(J)=J\cap[-N,N]^2.
\]
We say that $J$ has polynomial density if there exist
$\alpha\geq 1$, $\beta\geq 0$, and $\gamma_1,\gamma_2>0$ such that, for all
sufficiently large integers $N>1$,
\[
\frac{\gamma_1 N^{2/\alpha}}{(\log N)^\beta}
\leq \#Q_N(J)
\leq
\frac{\gamma_2 N^{2/\alpha}}{(\log N)^\beta}.
\]
In this case, we say that $J$ has polynomial density with exponent $\alpha$.
\begin{thm}[Fractal transference principle]\label{cor-FS}
Let $\Phi=\{\phi_i\colon i\in I\}$ be a $2$-decaying conformal IFS on $\Delta$. 
Let $S\subset I$. Suppose that $I$ has polynomial density with exponent $1$.
\begin{itemize}
\item[(a)] If $\overline{d}(S)>0$, then there exists $E_{S}(\Phi)\subset E(\Phi)$ 
such that \begin{equation}\label{equation1}\dim_{\rm H}E_{S}(\Phi)=\dim_{\rm H}E(\Phi)=1\ \text{ and }\ \overline{d}\left(\bigcup_{n\in\mathbb N}\bigcap_{x\in E_{S}(\Phi)}\{i_n(x)\}\cap S\right)=\overline{d}(S).\end{equation} 
In particular, 
\begin{equation}\label{equation-0}\dim_{\rm H}\left\{x\in E(\Phi)\colon \overline{d}(\{i_n(x)\colon n\in\mathbb N\}\cap S)=\overline{d}(S)\right\}=\dim_{\rm H}E(\Phi).\end{equation}
\item[(b)] If 
$d_{\rm B}(S)>0$, then
there exists $F_{S}(\Phi)\subset E(\Phi)$ 
such that \begin{equation}\label{equation1'}\dim_{\rm H}F_{S}(\Phi)=\dim_{\rm H}E(\Phi)=1\
\text{ and }\
d_{\rm B}\left(\bigcup_{n\in\mathbb N}\bigcap_{x\in F_{S}(\Phi)}\{i_n(x)\}\cap S\right)=d_{\rm B}(S).\end{equation}
In particular, 
\begin{equation}\label{equation-1}\dim_{\rm H}\left\{x\in E(\Phi)\colon d_{\rm B}(\{i_n(x)\colon n\in\mathbb N\}\cap S)=d_{\rm B}(S)\right\}=\dim_{\rm H}E(\Phi).\end{equation}
\end{itemize}\end{thm}

We next recall two fundamental results in multidimensional density combinatorics
that will be used in conjunction with our fractal transference principle.
The first is the following theorem of Furstenberg and Katznelson.
\begin{thm}[{\cite[Theorem~B]{FK78}}]\label{FK-thm} Let $k\in\mathbb N$ and suppose $S\subset\mathbb Z^k$ has positive upper Banach density. Then $S$ contains a homothetic copy of any finite subset of $\mathbb Z^k$.
\end{thm}

Bergelson and Leibman \cite{BL} proved a vast generalization of Furstenberg-Katznelson's theorem, known as polynomial Szemer\'edi's theorem.
Let $\mathbb Q_0[X]$ denote the set of polynomials with rational coefficients taking integer values on the integers and vanishing at $X=0$. 
\begin{thm}[{\cite[Theorem~${\rm B}_0$]{BL}}]\label{thm-BL}
Let $S\subset\mathbb Z^k$,
$k\in\mathbb N$ be a set of positive upper Banach density and let 
$h_1,\ldots,h_\ell$, $\ell\geq2$ be elements of $\mathbb Q_0[X]$. For any $v_1,\ldots,v_\ell \in\mathbb Z^k$
there exist $n\in\mathbb N$ 
and $u\in\mathbb Z^k$ such that 
 \[\{u+h_j(n)v_j\colon 1\leq j\leq \ell\}\subset S.\]
 \end{thm}
Taking $h_1=\cdots=h_\ell=X$ recovers Furstenberg-Katznelson's theorem.
 Transferring Theorem~\ref{thm-BL} 
via Theorem~\ref{cor-FS} we obtain the following statement. 

\begin{thm}\label{cor-N}
Let $\Phi=\{\phi_i\colon i\in I\}$ be a $2$-decaying conformal IFS on $\Delta$. 
Let $S\subset I$ have positive upper Banach density. Suppose that $I$ has polynomial density with exponent $1$.
  For any $h_1,\ldots,h_\ell\in\mathbb Q_0[X]$,
  $v_1,\ldots,v_\ell\in \mathbb Z^2$, $\ell\geq2$  
  there exist $n\in\mathbb N$, $u\in\mathbb Z^2$ and
  $n_1,\ldots,n_\ell\in \mathbb N$
  such that  
  \[\{u+h_j(n)v_j\colon 1\leq j\leq \ell\}\subset S,\]
  and 
  for any $x\in F_S(\Phi)$, \[i_{n_j}(x)=u+h_j(n)v_j\ \text{ for }j=1,\ldots,\ell.\]
\end{thm}

\subsection{A relativized version}
Theorem~\ref{cor-N} does not apply to subsets of $\mathbb Z^2$ whose upper Banach density is zero.
To treat such sets, notably the prime-coordinate set $\mathbb P^2$,
 we establish a relative fractal transference principle.
Given two infinite sets $A\subset S\subset\mathbb Z^2$, the number 
\[\overline{d}(A|S)=\limsup_{N\to\infty}\frac{\# Q_N(A) }{\# Q_N(S)}\]
is called
{\it the upper density of $A$ relative to $S$}.
Our second main result is 
a relativized version of Theorem~\ref{cor-FS}.

\begin{thm}[Relative fractal transference principle]\label{cor-FS-R}
Let $\Phi=\{\phi_i\colon i\in I\}$ be a $2$-decaying conformal IFS on $\Delta$. 
Let $S\subset I$ have polynomial density with exponent $\alpha\geq1$.
For any $A\subset S$ with $\overline{d}(A|S)>0$, 
there exists a subset $E_{S,A}(\Phi)$ of $\{x\in E(\Phi) \colon \{i_n(x)\colon n\in\mathbb N\}\subset  S \}$ such that \begin{equation}\label{equation4}\dim_{\rm H}E_{S,A}(\Phi)=\dim_{\rm H}\{x\in E(\Phi)\colon \{i_n(x)\colon n\in\mathbb N\}\subset S\}=\frac{1}{\alpha}\end{equation} and
\begin{equation}\label{equation5}\overline{d}\left(\bigcup_{n\in\mathbb N}\bigcap_{x\in E_{S,A}(\Phi)}\{i_n(x)\}\cap A|S\right)=\overline{d}(A|S).\end{equation}
In particular, 
\begin{equation}\label{equation6}\dim_{\rm H}\left\{
\begin{tabular}{l}
\!\!\!$x\in E(\Phi)\colon$\\
\!\!\!$\{i_n(x)\colon n\in\mathbb N\}\subset S$,\!\!\!\\
\!\!\!$\overline{d}(\{i_n(x)\colon n\in\mathbb N\}\cap A |S)$\!\!\!\\
\!\!\!$=\overline{d}(A|S)$\!\!\!
\end{tabular}
\right\}=\dim_{\rm H}
\left\{\begin{tabular}{l}
\!\!\!$x\in E(\Phi)\colon$\\
\!\!\!$\{i_n(x)\colon n\in\mathbb N\}\subset S$\!\!\!\end{tabular}
\right\}.\end{equation}
\end{thm}

 Green and Tao \cite{GT08} proved that the set $\mathbb P$ of primes contains arbitrarily long arithmetic progressions. Multidimensional extensions of Green-Tao's theorem were obtained by Cook, Magyar and  Titichetrakun \cite{CMT}, Fox and Zhao \cite{FZ15}, Tao and Ziegler \cite{TZ15}. Under the natural identification of $\mathbb C$ with $\mathbb R^2$,
the set $\mathbb P^2$ is regarded as the set of lattice points
whose two coordinates are prime.
\begin{thm}[{\cite[Theorem 1.4]{TZ15}}]
\label{TZ-thm} 
If $A\subset\mathbb P^2$ and $\overline{d}(A|\mathbb P^2)>0$,
then $A$ contains a homothetic copy of any finite subset of $\mathbb Z^2.$
\end{thm}
Since deleting finitely many points from $A$ does not change
$\overline d(A\mid\mathbb P^2)$, repeated applications of
Theorem~\ref{TZ-thm} yield infinitely many distinct homothetic copies.
By the prime number theorem, $\mathbb P^2$ has polynomial density with $\alpha=1$, $\beta=2$.  
Transferring Theorem~\ref{TZ-thm} 
via Theorem~\ref{cor-FS-R}  
we obtain the following statement.

\begin{thm}\label{GT-trans}
Let $\Phi=\{\phi_i:i\in I\}$ be a $2$-decaying conformal IFS on
$\Delta$, and suppose that $\mathbb P^2\subset I$.
For any $A\subset\mathbb P^2$ with
$
\overline d(A\mid\mathbb P^2)>0,$
there exists a set $ E_{\mathbb P^2,A}\subset E(\Phi)$ with Hausdorff dimension $1$ such that 
for any $x\in E_{\mathbb P^2,A}$, 

\begin{itemize}
\item $\{i_n(x)\colon n\in\mathbb N\}$ is contained in $\mathbb P^2$, and

\item 
$\{i_n(x)\colon n\in\mathbb N\}\cap A$
contains
infinitely many homothetic copies of any finite subset of $\mathbb Z^2$.
\end{itemize}
\end{thm}

\subsection{Outline of the proofs and organization of the paper}
\label{outlinepf}
Our proof is inspired by the  
surgery-and-elimination scheme developed in \cite{NT}, but the mechanism for preserving
Hausdorff dimension is different.  In~\cite{NT}, prescribed digits
are inserted into a one-dimensional continued-fraction construction
and the resulting set is compared with the original seed set through
an almost Lipschitz elimination map.  Here we instead establish a
direct dimension estimate for Moran-type constructions with sparse
deterministic insertions.  The proof uses conformality, the open set
condition, and planar packing estimates, and therefore does not rely
on the ordering properties specific to regular continued fractions.
This allows the argument to apply uniformly to a broad class of
conformal iterated function systems.


A different extension of \cite{NT} was obtained by
Guo et al.~\cite{GOQZ} for coordinatewise regular
continued fractions on $(0,1)^d$, $d\geq2$.  Their proof retains the
scheme of~\cite{NT}, but their product setting
introduces an anisotropic cylinder geometry.
Our systems treated in this paper are conformal, and hence isotropic and
 do not have a product structure.

The proofs of Theorems~\ref{cor-FS} and~\ref{cor-FS-R} follow
a common construction.
We first choose successive digits from suitable disjoint finite
sets, and then insert prescribed finite blocks at the same
positions in every sequence, without repeating digits.
These blocks realize the required density properties.
Their positions are chosen so that the cumulative sum of the
logarithms of the norms of inserted digits is negligible compared
with the corresponding sum for the original digits.

Proposition~\ref{prop:direct-deterministic-insertions} shows that
such insertions retain the required Hausdorff-dimension lower bound.
Its proof controls how many suitably selected images of $\Delta$
under finite compositions meet a small ball.
The open set condition and estimates for conformal maps provide
disjoint balls inside these images.
Comparing their areas gives a bound on the number of images,
which yields the measure estimates needed for the dimension bound.

For Theorem~\ref{cor-FS-R}, the original digits are chosen from
a subset of $S$ of zero relative upper density that still provides
enough choices for the dimension estimate.
We then insert the digits of $A$ outside this subset.
For Theorem~\ref{cor-FS}(b), we instead select widely separated
squares on which the density of $S$ approaches its upper Banach
density, while reserving enough digits for the initial construction.
The proof of Theorem~\ref{cor-FS}(a) uses centered squares.

Section~2 provides preliminary estimates and the dimension upper
bound.
Section~3 establishes the insertion estimate, and Section~4
completes the proofs of Theorems~\ref{cor-FS} and~\ref{cor-FS-R}.
Section~5 presents applications to the Hurwitz and
Kaneiwa--Shiokawa--Tamura continued fractions and the conformal
IFS of Mauldin and Urba\'nski.

\section{Geometric preliminaries}
 In Section \ref{univ-sec} we summarize basic properties of univalent functions and conformal IFSs
and in Section~2.2 we recall an upper bound for
the Hausdorff dimension of digit-restricted limit sets.

\subsection{Properties of univalent functions and conformal IFSs}\label{univ-sec}

The next lemma is an elementary consequence of Koebe’s distortion theorem; see, for example, \cite[Theorem~1.4]{CG93}.
\begin{lem}[{\cite[Lemma~2.2]{NT2}}]\label{distortion-lem}Let $\Omega\subset\mathbb C$ be a region and let $A$ be a compact subset of $\Omega$.
  There exists a constant $K\geq1$ such that for every univalent function $\phi\colon\Omega\to\mathbb C$ and all $x_1$, $x_2\in A$ we have
\[\frac{|D\phi(x_1)|}{|D\phi(x_2)|}\leq K.\]\end{lem}
For $x\in\mathbb C$ and $\delta>0$ let
$B(x, \delta)=\{y\in\mathbb C\colon |x-y|< \delta\}$.
\begin{lem}[{\cite[Lemma~2.3]{NT2}}]\label{conf-lem}
Let $K\geq 1$,
let $\Omega\subset \mathbb C$ be a region and let $\phi\colon \Omega\to \mathbb C$ be univalent. 
For any $x\in \Omega$ and $\delta>0$ such that \[ B(x, \delta)\subset \Omega \ \text{ and }
\ \sup_{x_1,x_2\in  B(x, \delta)}\frac{|D\phi(x_1)|}{|D\phi(x_2)|}\leq K,\]
we have
\[ B\left(\phi(x), \frac{\delta}{3K}|D\phi(x)| \right)\subset \phi( B(x, \delta) ).\]\end{lem}

For a set $A\subset\mathbb C$
let ${\rm diam}(A)=\sup\{|x-y|\colon x,y\in A\}$.
\begin{lem}
\label{proper2-new}
Let $\{\phi_i\colon i\in I\}$ be a 
$2$-decaying conformal IFS on $\Delta$.
Let $\Delta'\subset\mathbb C$ be a convex open set such that $\Delta\subset\Delta'$ and $\overline{\Delta'}\subset\widetilde\Delta$.
For any $\theta\in(0,1)$
there exists $C\geq1$ such that 
for all $n\in\mathbb N$ and $(i_1,\ldots, i_n)\in  I^n$ such that all the indices $i_1,\ldots,i_n$ are distinct, we have
\[C^{- 1}\prod_{k=1}^n\frac{1}{\|i_k\|_\infty^{2+\theta}}\le {\rm diam}(\phi_{i_1\cdots i_n}(\Delta'))\le C\prod_{k=1}^n\frac{1}{\|i_k\|_\infty^{2-\theta}}.\]
\end{lem}
\begin{proof}Fix $x_0\in {\rm int}\Delta$ and $\delta_0>0$ such that $B(x_0,\delta_0)\subset {\rm int}\Delta$.  
By Lemma~\ref{distortion-lem} and the convexity of $\Delta'$, there is  $K\geq1$ such that for any $n\in\mathbb N$ and any $(i_1,\ldots, i_n)\in  I^n$,
\begin{equation}\label{diam-eq1}{\rm diam}(\phi_{i_1\cdots i_n}(\Delta'))\leq K\cdot{\rm diam}(\Delta')|D\phi_{i_1\cdots i_n}(x_0)|.\end{equation}
By Lemma~\ref{conf-lem} we have 
\[ B\left(\phi_{i_1\cdots i_n}(x_0), \frac{\delta_0}{3K}|D\phi_{i_1\cdots i_n}(x_0)| \right)\subset \phi_{i_1\cdots i_n}( B(x_0, \delta_0) ),\] 
and in particular
\begin{equation}\label{diam-eq2}{\rm diam}(\phi_{i_1\cdots i_n}(\Delta'))\geq \frac{\delta_0}{3K}|D\phi_{i_1\cdots i_n}(x_0)|.\end{equation}

 Since $\{\phi_i\colon i\in I\}$ is 
$2$-decaying, there is $N\geq1$ such that
for any $i\in I$ with $\|i\|_\infty>N$ we have
\begin{equation}\label{diam-eq3}\frac{1}{\|i\|_\infty^{2+\theta}}\leq\frac{C_0^{-1} }{\|i\|_\infty^{2} }\leq \frac{C_0 }{\|i\|_\infty^{2} }\leq\frac{1}{\|i\|_\infty^{2-\theta}}.\end{equation}
Combining \eqref{diam-eq1} and the last inequality in \eqref{diam-eq3} yields
\[ {\rm diam}(\phi_{i_1\cdots i_n}(\Delta'))\le K\cdot{\rm diam}(\Delta')C_0^{\#\{1\leq k\leq n\colon \|i_k\|_\infty\leq N\}}\prod_{k=1}^n\frac{1}{\|i_k\|_\infty^{2-\theta}}.\]Combining \eqref{diam-eq2} and the first inequality in \eqref{diam-eq3} yields
\[ {\rm diam}(\phi_{i_1\cdots i_n}(\Delta'))\geq
\frac{\delta_0}{3K}C_0^{-\#\{1\leq k\leq n\colon \|i_k\|_\infty\leq N\}}\prod_{k=1}^n\frac{1}{\|i_k\|_\infty^{2+\theta}}.\]
Since all $i_1,\ldots,i_n$ are distinct, we have 
$\#\{1\leq k\leq n\colon \|i_k\|_\infty\leq N\}\leq (2N+1)^2$.
Hence the desired assertion follows.
\end{proof}

\subsection{Upper bound of Hausdorff dimension}\label{conv-sec}
Let $\Phi=\{\phi_i\colon i\in I\}$ be an infinite conformal IFS.
Fractal sets we are going to consider are contained in sets of the form
\[G_\Phi(S)=\left\{z\in L'(\Phi)\colon i_n(z)\in S\text{ for all $n\in\mathbb N$}\text{ and }\lim_{n\to\infty}\|i_n(z)\|_{\infty}=\infty\right\},\]
where $S$ is an infinite subset of $I$. We give an upper bound on the Hausdorff dimension of these sets.
For an infinite subset $S$ of $\mathbb Z^2$, define the convergence exponent of $S$ by
\[\tau(S)=\inf\left\{t\geq0\colon\sum_{i\in S} \|i\|_{\infty}^{-t}<\infty\right\}.\]







\begin{pro}\label{upper-d}  
Let $\Phi=\{\phi_i\colon i\in I\}$ be a 
$2$-decaying conformal IFS. For any infinite subset $S$ of $I$ we have \[\dim_{\rm H}G_\Phi(S)\leq\frac{\tau(S)}{2}.\]
If moreover $S$ has polynomial density with exponent $\alpha\geq1$, then 
\[\dim_{\rm H}G_\Phi(S)\leq\frac{1 }{\alpha}.\]\end{pro}

\begin{proof} 
The first inequality follows from the proof of
\cite[Proposition~3.1]{NT2}. Although that proposition was stated
under a stronger boundary condition, its proof uses only
(A1)--(A3) and the derivative estimate \eqref{CF-der}, and therefore
applies under our assumptions.
To prove the second inequality, let $s>2/\alpha$.
Since $S$ has polynomial density with exponent
$\alpha$, there exists $C>0$ such that
\[
\#Q_N(S)\leq CN^{2/\alpha}
\]
for all sufficiently large $N$.
Then we have
\[
\begin{aligned}
\sum_{\substack{i\in S\\
2^k\leq\|i\|_\infty<2^{k+1}}}
\|i\|_{\infty}^{-s}
&\leq 2^{-ks}\#Q_{2^{k+1}}(S)\\
&\leq C2^{2/\alpha}
2^{-k(s-2/\alpha)}
\end{aligned}
\]
for all sufficiently large $k$.
Since $s>2/\alpha$, it follows that
$
\sum_{i\in S}\|i\|_{\infty}^{-s}<\infty,
$
which implies $\tau(S)\leq s$. Letting
$s\to 2/\alpha$, we obtain
$\tau(S)\leq 2/\alpha$, as required.\end{proof}

\subsection{Negligible overlaps}
We record the relation between first-level overlaps and points
with more than one address.

\begin{lem}\label{lem:negligible-overlaps}
Let $\Phi$ be an infinite conformal IFS satisfying (A1)--(A3).
Then
\[
L(\Phi)\setminus L'(\Phi)
=\bigcup_{\omega\in I^{<\mathbb N}}\phi_\omega(\mathcal O(\Phi)),
\]
where $I^{<\mathbb N}$ denotes the set of all finite words over $I$, including
the empty word. In particular, (A4) is equivalent to
\[
\dim_{\rm H}(L(\Phi)\setminus L'(\Phi))=0.
\]

\end{lem}
\begin{proof}
Let $x\in L(\Phi)\setminus L'(\Phi)$, and choose two distinct
addresses of $x$. If $\omega$ is their maximal common initial word,
then $x=\phi_\omega(y)$ for some
$y\in\phi_i(L(\Phi))\cap\phi_j(L(\Phi))$ with $i\neq j$.
Conversely, any point in $\phi_\omega(\mathcal O(\Phi))$
has at least two addresses. This proves the asserted equality.
Since $I^{<\mathbb N}$ is countable and each $\phi_\omega$ is Lipschitz on
$\Delta$, (A4) implies that the set of non-unique addresses
has Hausdorff dimension zero. The converse follows from
$\mathcal O(\Phi)\subset L(\Phi)\setminus L'(\Phi)$.
\end{proof}

\section{Dimension estimates under sparse deterministic insertions}
\label{sec:dimension-insertions}

In this section we develop a general mechanism for estimating the
Hausdorff dimension of Moran-type sets obtained by inserting prescribed
finite blocks into symbolic sequences.
Moran-type constructions for non-autonomous conformal iterated function
systems have been developed in considerable generality; see, for example,
Rempe-Gillen and Urba\'nski~\cite{RU}.
Here we require a more specific estimate adapted to sparse deterministic
insertions.

The construction consists of two parts.  First, at the $n$-th level we
choose one symbol from a finite alphabet $I^{(n)}$ contained in an
annulus of radius comparable to $t^n$.  These symbols form a
\emph{seed sequence}.  We then insert a prescribed finite word $V_n$
after the $n$-th seed symbol.

The main point of this section is that the inserted blocks do not
affect the Hausdorff-dimension lower bound as long as their cumulative
logarithmic contraction cost is negligible compared with the
quadratic scale $n^2$ of the seed construction.  No density assumption
is imposed in this section.  Density-theoretic constructions will be
carried out in Section~4.

\subsection{Seed Moran constructions}

Let
$\Phi=\{\phi_i:i\in I\}$
be a $2$-decaying conformal IFS on $\Delta$.
Fix real numbers
$t>L>2.$
For each $n\geq1$, let $I^{(n)}$ be a non-empty finite subset of $I$
such that
\begin{equation}
 t^n\leq \|i\|_\infty\leq Lt^n
 \qquad (i\in I^{(n)}).
 \label{eq:integrated-annular-alphabets}
\end{equation}
Put
\[
 \mathbf I^\infty
 :=
 \prod_{n=1}^\infty I^{(n)}.
\]
The corresponding seed Moran set is defined by
$ R(\Phi,\mathbf I)
 :=
 \Pi(\mathbf I^\infty).$
Set
\begin{equation}
 D_*
 :=
 \liminf_{n\to\infty}
 \frac{\displaystyle\sum_{j=1}^n\log \#I^{(j)}}
      {n(n+1)\log t}.
 \label{eq:integrated-seed-dimension}
\end{equation}
Thus $D_*$ measures the asymptotic growth of the number of choices in
the seed construction relative to its geometric contraction scale.

\subsection{Deterministic insertion schemes}

For each $n\geq1$, let
\[
 V_n=v_{n,1}\cdots v_{n,m_n}
\]
be a finite word over $I$, where the empty word is allowed.
Assume that all symbols appearing in the words $(V_n)_{n\geq1}$
are pairwise distinct and that none of them belongs to
$\bigcup_{n\geq1}I^{(n)}.$
Define the cumulative logarithmic cost of the inserted symbols by
\begin{equation}
 \Lambda_n
 :=
 \sum_{j=1}^n\sum_{\ell=1}^{m_j}
 \log\|v_{j,\ell}\|_\infty .
 \label{eq:integrated-insertion-cost}
\end{equation}
Our basic negligibility assumption is
\begin{equation}
\lim_{n\to \infty} \frac{\Lambda_n}{n^2}=0.
 \label{eq:integrated-negligible-cost}
\end{equation}
For a finite seed word
$
 \omega=\omega_1\cdots\omega_n
 \in \prod_{j=1}^n I^{(j)},$
define its expanded word by
\begin{equation}
 \widehat{\omega}
 :=
 \omega_1V_1\omega_2V_2\cdots\omega_nV_n.
 \label{eq:integrated-expanded-word}
\end{equation}
For
$
 \omega=(\omega_n)_{n\geq1}\in\mathbf I^\infty,$
define the insertion operator
\[
 \mathcal I(\omega)
 :=
 \omega_1V_1\omega_2V_2\cdots.
\]
Put
$\widehat{\Pi}:=\Pi\circ\mathcal I$
and
\begin{equation}
 \widehat R=\widehat R(\Phi,\mathbf I)
 :=
 \widehat{\Pi}(\mathbf I^\infty).
 \label{eq:integrated-target-set}
\end{equation}

The following proposition is the main result of this section.  It
shows that insertions satisfying
\eqref{eq:integrated-negligible-cost} preserve the dimension supplied
by the seed alphabets.

\begin{pro}
\label{prop:direct-deterministic-insertions}
Under the assumptions above,
\begin{equation}
 \dim_{\rm H}\bigl(\widehat R\cap L'(\Phi)\bigr)
 \geq D_*,
 \qquad
 \widehat R\cap L'(\Phi)\subset E(\Phi).
 \label{eq:integrated-main-conclusion}
\end{equation}
\end{pro}

\subsection{Proof of Proposition~\ref{prop:direct-deterministic-insertions}}

For a finite word $\tau$ occurring in $\mathbf I^\infty$, let
\[
 [\tau]
 :=
 \{\omega\in\mathbf I^\infty:
   \omega|_{|\tau|}=\tau\}.
\]
Choose a convex open set $\Delta'$ such that
$\Delta\subset\Delta' $ and $\overline{\Delta'}\subset  \widetilde{\Delta},$
and fix
$ x_*\in\operatorname{int}\Delta.$
Since $t>L$, the annuli in
\eqref{eq:integrated-annular-alphabets} are pairwise disjoint.
Hence, the seed symbols in every admissible finite word are pairwise
distinct.  Together with the assumptions on the inserted symbols,
this implies that all symbols occurring in each expanded word
$\widehat{\omega}$ are pairwise distinct.  Therefore Lemma~\ref{proper2-new} is
applicable to $\phi_{\widehat{\omega}}$.


\begin{lem}
\label{lem:integrated-gibbs}
There exists a Borel probability measure $\mu$ on
$\mathbf I^\infty$ with the following property:
for every $0<s<D_*$ there exists $C_s>0$ such that, for every
$n\geq1$ and every
$\omega\in\prod_{j=1}^n I^{(j)},$
we have
\begin{equation}
 \mu([\omega])
 \leq
 C_s
 \bigl(
 \diam(\phi_{\widehat{\omega}}(\Delta'))
 \bigr)^s.
 \label{eq:integrated-cylinder-mass}
\end{equation}
\end{lem}

\begin{proof}
For each $n\geq1$, give $I^{(n)}$ the uniform probability
distribution.  By Kolmogorov's extension theorem, there exists a
Borel probability measure $\mu$ on $\mathbf I^\infty$ satisfying
\begin{equation}
 \mu([\omega_1,\ldots,\omega_n])
 =
 \prod_{j=1}^n\frac{1}{\#I^{(j)}}.
 \label{eq:integrated-product-measure}
\end{equation}
Fix $0<s<D_*$.  Choose $\theta\in(0,1)$ and $\eta>0$ sufficiently
small that
\begin{equation}
 s\left(1+\frac{\theta}{2}+\eta\right)<D_*.
 \label{eq:integrated-parameter-choice}
\end{equation}
By the lower estimate in Lemma~\ref{proper2-new},
\eqref{eq:integrated-annular-alphabets}, and
\eqref{eq:integrated-insertion-cost}, for $|\omega|=n$ we have
\begin{align}
 \diam(\phi_{\widehat{\omega}}(\Delta'))
 &\geq
 C_\theta^{-1}
 \prod_{j=1}^n
 \frac{1}{(Lt^j)^{2+\theta}}
 \prod_{j=1}^n
 \prod_{\ell=1}^{m_j}
 \frac{1}{\|v_{j,\ell}\|_\infty^{2+\theta}}
 \notag\\
 &=
 C_\theta^{-1}
 L^{-(2+\theta)n}
 t^{-n(n+1)(1+\theta/2)}
 \exp\bigl(-(2+\theta)\Lambda_n\bigr).
 \label{eq:integrated-diameter-lower-pre}
\end{align}
By \eqref{eq:integrated-negligible-cost},
after increasing the constant to deal with finitely many initial
levels, we obtain
\begin{equation}
 C_\theta^{-1}
 L^{-(2+\theta)n}
 \exp\bigl(-(2+\theta)\Lambda_n\bigr)
 \geq
 C_{\theta,\eta}^{-1}
 t^{-\eta n(n+1)}.
 \label{eq:integrated-cost-absorbed}
\end{equation}
Consequently,
\begin{equation}
 \diam(\phi_{\widehat{\omega}}(\Delta'))
 \geq
 C_{\theta,\eta}^{-1}
 t^{-n(n+1)(1+\theta/2+\eta)}.
 \label{eq:integrated-diameter-lower}
\end{equation}
Put
$
 A_n:=\sum_{j=1}^n\log\#I^{(j)}.
$
By \eqref{eq:integrated-seed-dimension} and
\eqref{eq:integrated-parameter-choice}, there exist $\delta>0$ and
$n_0\geq1$ such that for any $n\geq n_0,$
\begin{equation}
 A_n
 \geq
 \left(
 s\left(1+\frac{\theta}{2}+\eta\right)+\delta
 \right)n(n+1)\log t.
 \label{eq:integrated-entropy-lower}
\end{equation}
Combining
\eqref{eq:integrated-product-measure},
\eqref{eq:integrated-diameter-lower}, and
\eqref{eq:integrated-entropy-lower} yields for $n\geq n_0$,
\begin{align}
 \frac{\mu([\omega])}
 {\diam(\phi_{\widehat{\omega}}(\Delta'))^s}
 &\leq
 C_{\theta,\eta}
 \exp(-A_n)
 t^{s(1+\theta/2+\eta)n(n+1)}
 \notag\\
 &\leq
 C_{\theta,\eta}
 t^{-\delta n(n+1)}.
 \label{eq:integrated-mass-ratio}
\end{align}
The finitely many levels $n<n_0$ can be absorbed into the constant.
This proves \eqref{eq:integrated-cylinder-mass}.
\end{proof}
For a finite admissible seed word $\omega$, put
\begin{equation}
 \rho(\omega)
 :=
 \bigl|D\phi_{\widehat{\omega}}(x_*)\bigr|.
 \label{eq:integrated-domega}
\end{equation}
Set
\[
 r_0
 :=
 \min_{i\in I^{(1)}}
 \bigl|D\phi_{iV_1}(x_*)\bigr|
 >0.
\]

For $x\in\widehat R$ and $r\in(0,r_0)$, let
$\mathcal V(x,r)$ be the set of all finite admissible seed words
$\omega$ satisfying

\begin{enumerate}
\item
$
 B(x,r)\cap\phi_{\widehat{\omega}}(\Delta)\neq\emptyset;$

\item
$
 \rho(\omega)\geq r,\ \text{and}\ \text{there exists}\ 
 a\in I^{(|\omega|+1)} \ \text{such that}\ 
 \rho(\omega a)<r.
$

\end{enumerate}
Here,
$
 \widehat{\omega a}
 =
 \widehat{\omega}\,aV_{|\omega|+1}.$
Since each $I^{(n)}$ is finite,
\begin{equation}
 \min\{|\omega|:\omega\in\mathcal V(x,r)\}
 \longrightarrow\infty
 \qquad(r\to0).
 \label{eq:integrated-min-length}
\end{equation}
Let $\mathcal V'(x,r)$ be the set of
$\omega\in\mathcal V(x,r)$ for which there exists
$\tau\in\mathcal V(x,r)$ such that
$
 |\tau|<|\omega|
 $ and $
 [\omega]\subset[\tau].
$
Set
$
 \mathcal V^*(x,r)
 :=
 \mathcal V(x,r)\setminus\mathcal V'(x,r).
$


\begin{lem}
\label{lem:integrated-adjust}
If $x\in\widehat R$, $r\in(0,r_0)$, and
$\omega\in\mathcal V'(x,r)$, then there exists
$\tau\in\mathcal V^*(x,r)$ such that
$
 [\omega]\subset[\tau]
$ and $
 |\tau|<|\omega|.
$
\end{lem}

\begin{proof}
Starting with $\omega\in\mathcal V'(x,r)$, choose a shorter word
$\tau_1\in\mathcal V(x,r)$ with
$[\omega]\subset[\tau_1].$
If $\tau_1\notin\mathcal V^*(x,r)$, choose a still shorter word with
the same property.  Since word lengths are non-negative integers,
the procedure terminates after finitely many steps at an element of
$\mathcal V^*(x,r)$.
\end{proof}


\begin{lem}
\label{lem:integrated-cover}
Let $x\in\widehat R$ and $r\in(0,r_0)$.  Then the following hold.

\begin{enumerate}
\item
For distinct $\omega,\tau\in\mathcal V^*(x,r)$,
\begin{equation}
 \phi_{\widehat{\omega}}(\operatorname{int}\Delta)
 \cap
 \phi_{\widehat{\tau}}(\operatorname{int}\Delta)
 =
\emptyset.
 \label{eq:integrated-open-set-disjointness}
\end{equation}

\item
\begin{equation}
 \widehat{\Pi}^{-1}(B(x,r))
 \subset
 \bigcup_{\tau\in\mathcal V^*(x,r)}[\tau].
 \label{eq:integrated-symbolic-cover}
\end{equation}
\end{enumerate}
\end{lem}

\begin{proof}
Distinct elements of $\mathcal V^*(x,r)$ are prefix-incomparable.
Suppose their first differing seed symbols occur at level $j$.
Their expanded words have the same prefix through level $j-1$ and
have different symbols at level $j$.  Hence, the open set condition
gives
\eqref{eq:integrated-open-set-disjointness}.

Let $\omega\in\mathbf I^\infty$ satisfy $\widehat{\Pi}(\omega)\in B(x,r).$
Since the derivatives of the expanded prefix maps tend to zero and
$r<r_0$, there exists $n\geq1$ such that
$
 \rho(\omega|_n)\geq r>\rho(\omega|_{n+1}).$
Moreover,
$
 \widehat{\Pi}(\omega)
 \in
 \phi_{\widehat{\omega|_n}}(\Delta),$
and hence
$
 \omega|_n\in\mathcal V(x,r).
$
If $\omega|_n\in\mathcal V^*(x,r)$ there is nothing to prove.
Otherwise Lemma~\ref{lem:integrated-adjust} gives
$\tau\in\mathcal V^*(x,r)$ such that
$[\omega|_n]\subset[\tau].$
In either case $\omega\in[\tau]$, proving
\eqref{eq:integrated-symbolic-cover}.
\end{proof}


\begin{lem}
\label{lem:integrated-stopping-estimates}
Fix $\theta\in(0,1)$ and $\eta>0$, and put
\begin{equation}
 a:=1+\frac{\theta}{2}+\eta,
 \qquad
 b:=1-\frac{\theta}{2},
 \qquad
 \kappa:=\frac{b}{a}.
 \label{eq:integrated-abkappa}
\end{equation}
Then there exist constants $C,c>0$ and
$r_1\in(0,r_0]$ such that, for every $x\in\widehat R$,
every $r\in(0,r_1)$, and every $\omega\in\mathcal V(x,r)$,
\begin{equation}
 \diam(\phi_{\widehat{\omega}}(\Delta'))
 \leq
 Cr^\kappa
 t^{c\sqrt{\log(1/r)}},
 \label{eq:integrated-stopping-diameter}
\end{equation}
and
\begin{equation}
 \#\mathcal V^*(x,r)
 \leq
 Cr^{-2(1-\kappa)}
 t^{c\sqrt{\log(1/r)}}.
 \label{eq:integrated-stopping-number}
\end{equation}
\end{lem}

\begin{proof}
By Lemma~\ref{distortion-lem}, there exists $K\geq1$ such that
\begin{equation}
 \sup_{z_1,z_2\in\Delta'}
 \frac{|D\phi_u(z_1)|}{|D\phi_u(z_2)|}
 \leq K
 \label{eq:integrated-bounded-distortion}
\end{equation}
for every finite word $u$ over $I$.
Choose $\delta_*>0$ such that
$B(x_*,\delta_*)\subset\operatorname{int}\Delta$.
By Lemma~\ref{conf-lem} and the convexity of $\Delta'$,
there exists $C_1\geq1$ such that
\begin{equation}
 C_1^{-1}\rho(\omega)
 \leq\diam(\phi_{\widehat\omega}(\Delta'))
 \leq C_1\rho(\omega)
 \label{eq:integrated-derivative-diameter}
\end{equation}
for every finite admissible seed word $\omega$.

Put $a=1+\theta/2+\eta$ and $b=1-\theta/2$, as in
\eqref{eq:integrated-abkappa}.
By Lemma~\ref{proper2-new} and
\eqref{eq:integrated-annular-alphabets}, there exists $C_2\geq1$
such that, for $n=|\omega|$,
\begin{equation}
 \diam(\phi_{\widehat\omega}(\Delta'))
 \leq C_2t^{-bn(n+1)}.
 \label{eq:integrated-diameter-upper}
\end{equation}
Here the factors corresponding to the inserted symbols are at most
one. On the other hand, the estimate
\eqref{eq:integrated-diameter-lower} yields a constant $C_3\geq1$
such that
\begin{equation}
 \diam(\phi_{\widehat\omega}(\Delta'))
 \geq C_3^{-1}t^{-an(n+1)}.
 \label{eq:integrated-diameter-two-sided-lower}
\end{equation}
The constants $C_1,C_2,C_3$ are independent of $\omega$.

Let $x\in\widehat R$, $r\in(0,r_0)$, and
$\omega\in\mathcal V(x,r)$. Write $n=|\omega|$.
By the definition of $\mathcal V(x,r)$,
\eqref{eq:integrated-derivative-diameter} and
\eqref{eq:integrated-diameter-upper}, there exists $C_4\geq1$
such that
$r\leq \rho(\omega)\leq C_4t^{-bn(n+1)}.$
Choose $r_1\in(0,r_0]$ so small that
$r_1<\min\{e^{-1},C_4^{-1}\}$.
Then $bn(n+1)\log t\leq 2\log(1/r)$ for $0<r<r_1$,
and hence
\begin{equation}
 n+1\leq C_5\sqrt{\log(1/r)}
 \label{eq:integrated-stopping-length}
\end{equation}
with a constant $C_5>0$ independent of $x,r,\omega$.

Choose $a_0\in I^{(n+1)}$ such that $\rho(\omega a_0)<r$.
Applying \eqref{eq:integrated-diameter-two-sided-lower}
to $\omega a_0$ and then using
\eqref{eq:integrated-derivative-diameter}, we have
\[
 C_3^{-1}t^{-a(n+1)(n+2)}
 \leq \diam(\phi_{\widehat{\omega a_0}}(\Delta'))
 \leq C_1\rho(\omega a_0)<C_1r.
\]
Hence, since $\kappa=b/a$, there exists $C_6>0$ such that
\[
 t^{-b(n+1)(n+2)}\leq C_6r^\kappa.
\]
Combining this with \eqref{eq:integrated-diameter-upper}
and \eqref{eq:integrated-stopping-length}, we obtain
\begin{align*}
 \diam(\phi_{\widehat\omega}(\Delta'))
 &\leq C_2t^{-bn(n+1)}\\
 &=C_2t^{2b(n+1)}t^{-b(n+1)(n+2)}\\
 &\leq C_7r^\kappa
 t^{c_1\sqrt{\log(1/r)}},
\end{align*}
where $C_7>0$ and $c_1>0$ are constants independent of
$x,r,\omega$. This proves
\eqref{eq:integrated-stopping-diameter} with $C_7,c_1$
in place of $C,c$.

Since $B(x,r)$ intersects $\phi_{\widehat\omega}(\Delta)$
for every $\omega\in\mathcal V(x,r)$, and since
$0<\kappa<1$, it follows that, after decreasing $r_1$
if necessary,
\begin{equation}
 \bigcup_{\omega\in\mathcal V(x,r)}
 \phi_{\widehat\omega}(\Delta)
 \subset B\left(x,C_8r^\kappa
 t^{c_1\sqrt{\log(1/r)}}\right)
 \label{eq:integrated-union-contained}
\end{equation}
for some constant $C_8>0$.

By Lemma~\ref{conf-lem} and
\eqref{eq:integrated-bounded-distortion}, there exists $c_*>0$
such that
\begin{equation}
 B\left(\phi_{\widehat\omega}(x_*),c_*\rho(\omega)\right)
 \subset\phi_{\widehat\omega}(\operatorname{int}\Delta)
 \label{eq:integrated-inner-ball}
\end{equation}
for every finite admissible seed word $\omega$.
For $\omega\in\mathcal V^*(x,r)$, this ball has radius
at least $c_*r$. These balls are pairwise disjoint by
Lemma~\ref{lem:integrated-cover}(1), and all are contained
in the ball on the right-hand side of
\eqref{eq:integrated-union-contained}.
Comparing planar Lebesgue measures, we obtain
\[
 \#\mathcal V^*(x,r)\,\pi(c_*r)^2
 \leq \pi C_8^2r^{2\kappa}
 t^{2c_1\sqrt{\log(1/r)}},
\]
and consequently
\[
 \#\mathcal V^*(x,r)
 \leq C_9r^{-2(1-\kappa)}
 t^{2c_1\sqrt{\log(1/r)}}
\]
for some $C_9>0$.
All constants introduced above depend only on the fixed data
and on $\theta,\eta$, and are independent of $x,r,\omega$.
Taking $C=\max\{C_7,C_9\}$ and $c=2c_1$,
we obtain both \eqref{eq:integrated-stopping-diameter}
and \eqref{eq:integrated-stopping-number}.
\end{proof}


\begin{proof}[Proof of Proposition~\ref{prop:direct-deterministic-insertions}]
If $D_*=0$, the dimension inequality in
\eqref{eq:integrated-main-conclusion} is immediate.  Assume therefore
that $D_*>0$.
Fix
$0<s<D_*$,
and choose $s_0$ such that
\begin{equation}
 s<s_0<D_*.
 \label{eq:integrated-s0-choice}
\end{equation}
Choose $\theta\in(0,1)$ and $\eta>0$ sufficiently small that
\begin{equation}
 s_0
 \left(
 1+\frac{\theta}{2}+\eta
 \right)
 <D_*
 \label{eq:integrated-s0-cylinder-condition}
\end{equation}
and, with $\kappa$ as in
\eqref{eq:integrated-abkappa},
\begin{equation}
 s_0\kappa-2(1-\kappa)>s.
 \label{eq:integrated-final-exponent-choice}
\end{equation}

Let $\mu$ be the probability measure from
Lemma~\ref{lem:integrated-gibbs}, applied with exponent $s_0$, and
define
$\nu
 :=
 \mu\circ\widehat{\Pi}^{-1}.$
Let $x\in\widehat R$ and let $r>0$ be sufficiently small.
By Lemma~\ref{lem:integrated-cover}(2),
Lemma~\ref{lem:integrated-gibbs}, and
Lemma~\ref{lem:integrated-stopping-estimates},
\begin{align}
 \nu(B(x,r))
 &\leq
 \sum_{\omega\in\mathcal V^*(x,r)}
 \mu([\omega])
 \notag\\
 &\leq
 C\,
 \#\mathcal V^*(x,r)
 \max_{\omega\in\mathcal V^*(x,r)}
 \diam(\phi_{\widehat{\omega}}(\Delta'))^{s_0}
 \notag\\
 &\leq
 C^{1+s_0}
 r^{s_0\kappa-2(1-\kappa)}
 t^{c(1+s_0)\sqrt{\log(1/r)}}.
 \label{eq:integrated-final-ball-estimate}
\end{align}
By
\eqref{eq:integrated-final-exponent-choice}, there exists
$\delta>0$ such that
$s_0\kappa-2(1-\kappa)
 \geq
 s+2\delta.$
Moreover,
\[
 t^{c(1+s_0)\sqrt{\log(1/r)}}
 \leq
 r^{-\delta}
\]
for all sufficiently small $r$.  Therefore
\begin{equation}
 \nu(B(x,r))
 \leq
 C^{1+s_0}r^{s+\delta}
 \leq
 C^{1+s_0}r^s.
 \label{eq:integrated-frostman-final}
\end{equation}
By the mass distribution principle,
$\dimH\widehat R\geq s.$
Since $0<s<D_*$ is arbitrary,
$
 \dimH\widehat R\geq D_*.
$
Condition (A4) implies
$ \dimH(L(\Phi)\setminus L'(\Phi))=0$
(see Lemma~\ref{lem:negligible-overlaps}).
Hence,
\[
 \dimH(\widehat R\cap L'(\Phi))
 \geq D_*.
\]

Finally, every expanded address constructed above has pairwise
distinct symbols.  Therefore, if $x\in\widehat R\cap L'(\Phi)$, its
unique address has pairwise distinct symbols.  Consequently,
\[
 \widehat R\cap L'(\Phi)\subset E(\Phi).
\]
This completes the proof.
\end{proof}

\subsection{A block criterion for negligible insertion cost}

For later applications it is convenient to give a simple sufficient
condition for
\eqref{eq:integrated-negligible-cost}.  It allows us to insert an
arbitrary finite block at a sufficiently sparse sequence of levels.

\begin{lem}
\label{lem:integrated-blockwise-cost}
Let
$0=M_0<M_1<M_2<\cdots,$
let $(\varepsilon_k)_{k\geq1}$ be a decreasing sequence of positive
numbers converging to $0$, and let $(W_k)_{k\geq1}$ be finite subsets
of $I$.  Define
\[
 V_n=
 \begin{cases}
 \text{a fixed enumeration of }W_k,
     & n=M_k,\\
\emptyset,
     & n\notin\{M_k:k\geq1\}.
 \end{cases}
\]
If
\begin{equation}
 \prod_{i\in W_k}\frac{1}{\|i\|_\infty^2}
 \geq
 \left(
 \prod_{j=M_{k-1}+1}^{M_k}\frac{1}{t^{2j}}
 \right)^{\varepsilon_k}
 \qquad(k\geq1),
 \label{eq:integrated-block-product-condition}
\end{equation}
then
$\lim_{n\to \infty}\Lambda_n/n^2=0.$
\end{lem}

\begin{proof}
Taking logarithms in
\eqref{eq:integrated-block-product-condition}, we obtain
\begin{equation}
 \sum_{i\in W_k}\log\|i\|_\infty
 \leq
 \varepsilon_k(\log t)
 \sum_{j=M_{k-1}+1}^{M_k}j.
 \label{eq:integrated-block-log-bound}
\end{equation}
Let $q=q(n)$ be determined by
$M_q\leq n<M_{q+1}.$
Fix $K\geq1$.  For all sufficiently large $n$, using the monotonicity
of $(\varepsilon_k)$ and
\eqref{eq:integrated-block-log-bound}, we have
\[
 \Lambda_n
 \leq
 \sum_{k=1}^K\sum_{i\in W_k}\log\|i\|_\infty
 +
 \varepsilon_K(\log t)
 \sum_{j=M_K+1}^{M_q}j.
\]
Dividing by \(n^2\) and letting \(n\to\infty\), we obtain
\[
 \limsup_{n\to\infty}
 \frac{\varepsilon_K(\log t)}{n^2}
 \sum_{j=M_K+1}^{M_q}j
 \leq
 \frac{\varepsilon_K\log t}{2}.
\]
Letting $K\to\infty$ gives
${\Lambda_n}/{n^2}\longrightarrow0,$ which is required.
\end{proof}

\section{Proof of the main results}
\subsection{Construction of a sparse seed set}\label{C-const}
 Let $\Phi=\{\phi_i\colon i\in I\}$ be a $2$-decaying conformal IFS on $\Delta$. Let $S\subset I$ have polynomial density with exponent $\alpha\geq1$.
Define a function $\lambda\colon [1,\infty)\to\mathbb N$ by \begin{equation}\label{mu-def}\lambda(x)=k\ \text{ for } x\in[k!,(k+1)!),\ k\in\mathbb N.\end{equation}
\begin{lem}[{\cite[Lemma~3.2]{NT}}]\label{simple}For any $t>e$ there exists $\sigma>1$ such that \[\lambda(t^n)\leq (n\log t)^{\sigma}\ \text{ for every }n\in\mathbb N.\]\end{lem}

The next lemma allows us to extract from an arbitrary infinite subset of $\mathbb Z^2$ a set of zero relative density with definite growth bounds.

\begin{lem}\label{find-C}Let
$S$ be an infinite subset of $\mathbb Z^2$ and enumerate $S$ without repetition as $S=\{i_n\in\mathbb Z^2\colon n\in\mathbb N\}$,  $\|i_1\|_\infty\leq \|i_{2}\|_\infty\leq\cdots.$
There exists a subset $S_*$ of $S$ such that 
  for all sufficiently large integer $T\geq1$,  \begin{equation}\label{seed-eq3}\frac{1}{2\lambda(T)}\leq \frac{\#\{1\leq n\leq T\colon i_n\in S_* \}}{T }\leq \frac{3}{\lambda(T)},\end{equation} 
  and \[\overline{d}(S_*|S)=0.\]
\end{lem}

\begin{proof}
Let \[P=\bigcup_{k=2}^\infty\{k!+ik\colon i=0,\ldots, k!-1\}\subset\mathbb N.\]
 In the proof of \cite[Lemma~3.2]{NT} it was shown that for all sufficiently large $T\geq1$,
\[\begin{split}\frac{T}{2\lambda(T)}\leq\#(P\cap[1,T])&
\leq\frac{3T}{\lambda(T)}.\end{split}\]
Set \[S_*=\{i_n\in S\colon n\in P\}.\]
For every $T\geq1$ we have
$\#\{1\leq n\leq T\colon i_n\in S_* \}=\#(P\cap[1,T]).$
This verifies \eqref{seed-eq3}.

For each integer $N\geq1$ put
$T_N=\max\{n\geq1\colon \|i_n\|_\infty\leq N\}$ for $N$ satisfying $Q_N(S)\neq \emptyset$.
We claim that \begin{equation}\label{claim1}\#Q_N(S)=T_N .\end{equation}
Since the nondecreasing sequence $(||i_n||_{\infty})_{n\ge 1}$ tends to infinity, we have $n\leq T_N\Rightarrow ||i_n||_{\infty}\leq N.$ By the maximality of  $T_N$, we have  $n> T_N\Rightarrow ||i_n||_{\infty}> N.$ Combining this we have $Q_N(S)=\{i_1,..., i_{T_N}\}$.
\noindent Then \eqref{claim1} implies 
\begin{equation}\label{claim2}\#Q_N(S_*) =\#\{1\leq n\leq
T_N\colon i_n\in S_* \}.\end{equation}
By \eqref{claim1}, \eqref{claim2} and the upper bound in \eqref{seed-eq3}, for all sufficiently large $N$ we have
\[\frac{\#Q_N(S_*) }{\#Q_N(S) }=\frac{\#\{1\leq n\leq
T_N\colon i_n\in S_* \}}{T_N }
\leq\frac{3}{\lambda(T_N )}.\]
Letting $N\to\infty$ yields $\overline{d}(S_*|S)=0$ as required. \end{proof}

Let $\alpha\geq1$, $\beta\geq0$, $\gamma_1$, $\gamma_2>0$ be the constants in the definition of polynomial density of $S$ with exponent $\alpha$. 
Set 
\begin{equation}\label{L}
L=(12\gamma_1^{-1}\gamma_22^\beta)^{\alpha/2}>1.\end{equation}
Let $t>\max\{L,e\}$ satisfy
 \begin{equation}\label{seed-eq4}\lambda(\#Q_t(S) )\geq2.\end{equation}
Let $S_*\subset  S$ be the  infinite subset introduced in Lemma~\ref{find-C}. 
For each $n\geq1$, let
\[
I^{(n)}=Q_{Lt^n}(S^*)\setminus Q_{t^{n}-1}(S^*)
.\]


We endow each $I^{(n)}$ with the discrete topology, and $I^{\infty}$ with the product topology.
For a finite word $\omega$, let $|\omega|$ denote its word length. For $\underline{\omega}=(\omega_1,\omega_2,\ldots)\in I^\infty$ and $n\in\mathbb N$, let $\underline{\omega}|_n=\omega_1\cdots\omega_n$.
\begin{lem}
\label{Cardi}
There exists $\rho>0$ such that for any sufficiently large $n\ge 1,$
\begin{equation*}\begin{split}\#I^{(n)}\geq\frac{t^{\frac{2n}{\alpha}}}{(3n\log t)^{\rho }}.\end{split}\end{equation*}
In particular, 
\[
\liminf_{n\to\infty}
 \frac{\displaystyle\sum_{j=1}^n\log \#I^{(j)}}
      {n(n+1)\log t}\ge \frac{1}{\alpha}.
      \]
\end{lem}
\begin{proof}
Let $\Phi=\{\phi_i\colon i\in I\}$ be a $2$-decaying conformal IFS on $\Delta$.
Let $S\subset I$ have polynomial density with exponent $\alpha\ge 1$ and write 
\begin{equation}\label{alignment}S=\{i_n\colon n\in\mathbb N\},\ \|i_1\|_\infty\leq\|i_2\|_\infty\leq\cdots.\end{equation}  


The ordering in \eqref{alignment} implies that for all $N\geq1$, 
\[\#Q_N(S_*) 
=\#\{1\leq n\leq \#Q_N(S) \colon i_n\in S_*\}.\]
By this and \eqref{seed-eq3} we have 
\begin{equation}\label{seed-eq8}\frac{\#Q_N(S) }{2\lambda(\#Q_N(S))}\leq\#Q_N(S_*)
\leq \frac{3\#Q_N(S) }{\lambda(\#Q_N(S) )}.\end{equation}

 Since $t>L$ and $S$ has polynomial density, if $t$ is sufficiently large then for every $n\in\mathbb N$ we have
\begin{equation}\label{seed-eq0}\gamma_1\gamma_2^{-1}2^{-\beta}L^{-\frac{2}{\alpha}}\leq\frac{\# Q_{t^n}(S) }{\# Q_{Lt^n}(S) }\leq \gamma_1^{-1}\gamma_22^{\beta}L^{-\frac{2}{\alpha}}.\end{equation}
For each $n\in\mathbb N$, fix $k=k(n)\in\mathbb N$ such that
$k!\leq\#Q_{t^n}(S) <(k+1)!.$
If $n$ is sufficiently large so that
$k+2\geq (\gamma_1\gamma_2^{-1}2^{-\beta}L^{-\frac{2}{\alpha}})^{-1}$, then
\[k!\leq\#Q_{t^n}(S)\leq\#Q_{Lt^n}(S)\leq\frac{\# Q_{t^n}(S) }{\gamma_1\gamma_2^{-1}2^{-\beta}L^{-\frac{2}{\alpha}}}<(k+2)!.\]
By the definition of $\lambda$ in \eqref{mu-def}, this estimate implies
 \begin{equation}\label{seed-eq2}0\leq\lambda(\#Q_{Lt^n}(S))-\lambda(\#Q_{t^n}(S) )\leq 1.\end{equation}
Since $t>e$, Lemma~\ref{simple} gives $\sigma>1$ such that
for every $n\in\mathbb N$, after increasing $\sigma$ if necessary
\begin{equation}\label{seed-eq5}\lambda(\#Q_{t^n}(S) )\leq\lambda((2t^{n}+1)^2)\leq \lambda(t^{3n})\leq \gamma_1(3n\log t)^{\sigma }.\end{equation}

 Set $\rho=\beta+\sigma$. 
For each $n\in\mathbb N$ put
\[I^{(n)}=Q_{Lt^{n}}(S_*)\setminus Q_{t^{n}-1}(S_*).\]
We evaluate $\#I^{(n)}$.
Using \eqref{seed-eq8}, \eqref{seed-eq4} and then \eqref{seed-eq0}, \eqref{seed-eq2}, 
 for each integer $n\geq 1$ we have
 \[\begin{split}\#I^{(n)}\geq&
 \#Q_{Lt^{n}}(S_*) 
 -\#Q_{t^{n}}(S_* )\\
\geq&\frac{\#Q_{Lt^n}(S) }
{2\lambda(\#Q_{Lt^n}(S)) }-
\frac{3\#Q_{t^n}(S)}{\lambda(\#Q_{t^n}(S))}\\
\geq&\frac{\gamma_1\gamma_2^{-1} 2^{-\beta}L^{\frac{2}{\alpha}}
\#Q_{t^n}(S)}{2(\lambda(\#Q_{t^n}(S))+1)}-\frac{3\#Q_{t^n}(S)}
{\lambda(\#Q_{t^n}(S))}.\end{split}\]
By \eqref{seed-eq4}, the last expression is bounded from below by
\[\#Q_{t^n}(S) 
\left(\frac{\gamma_1\gamma_2^{-1}2^{-\beta}
L^{\frac{2}{\alpha}} }{3
\lambda(\#Q_{t^n}(S))}-\frac{3}
{\lambda(\#Q_{t^n}(S) )}\right),\]
which equals $\#Q_{t^n}(S)/\lambda(\#Q_{t^n}(S))$ by \eqref{L}.
Using \eqref{seed-eq5} and the assumption that $S$ has polynomial density with exponent $\alpha$, we obtain 
\begin{equation*}\begin{split}\#I^{(n)}\geq\frac{\#Q_{t^n}(S) }{\lambda(\#Q_{t^n}(S) )}\geq \gamma_1\frac{t^{\frac{2n}{\alpha}}}{(n\log t)^{\beta}}\cdot\frac{1}{\gamma_1(3n\log t)^{\sigma }}
\geq \frac{t^{\frac{2n}{\alpha}}}{(3n\log t)^{\rho }},\end{split}\end{equation*} as required.
\end{proof}


\subsection{Density-preserving insertions and proof of Theorems~\ref{cor-FS-R}   and~\ref{cor-FS}(a)}\label{gen-sec}
Let $\Phi=\{\phi_i\colon i\in I\}$ be a $2$-decaying conformal IFS on $\Delta$.
Let $S\subset I$ have polynomial density with exponent $\alpha$. 
Let $A\subset S$
satisfy $\overline{d}(A|S)>0$. We have
\begin{equation}\label{dim-equality}\overline{d}((A\setminus S_*)|S)=\overline{d}(A|S).\end{equation}

Let $(\varepsilon_k)_{k=1}^\infty$ be a decreasing sequence of positive real
numbers converging to $0$. Choose a strictly increasing sequence
$(N_k)_{k=1}^\infty$ of positive integers such that
\begin{equation}\label{eq-sec4-1}
\lim_{k\to\infty}
\frac{\#Q_{N_k}(A\setminus S_{\ast})}{\#Q_{N_k}(S)}=\overline{d}(A\setminus S_{\ast}| S)=\overline{d}(A| S).
\end{equation}
Set $N_0:=0$, and for each $k\ge 1$ define
\[W_k:=Q_{N_k}(A\setminus S_{\ast})\setminus Q_{N_{k-1}}(A\setminus S_{\ast}).\]
Then each $W_k$ is a finite subset of $A\setminus S_{\ast}$, the family
$\{W_k\colon k\ge 1\}$ is pairwise disjoint, and
\begin{equation}\label{eq-sec4-2}
\bigcup_{j=1}^k W_j = Q_{N_k}(A\setminus S_{\ast})
\ \text{for all}\ k\ge 1.
\end{equation}

We now construct a strictly increasing sequence $(M_k)_{k=1}^\infty$ inductively.
Set $M_0:=0$. Suppose that $M_0,\dots,M_{k-1}$ have already been chosen. Since
\[\sum_{j=M_{k-1}+1}^{\ell}j\to \infty\  \text{as}\ \ell\to \infty,\]
we can choose $M_k>M_{k-1}$ such that
\begin{equation}
\label{eq-sec4-4}
\prod_{i\in W_{k}}\frac{1}{\|i\|_\infty^{2}}\geq\left(\prod_{j=M_{k-1}+1}^{M_k}\frac{1}{t^{2 j}} \right)^{\varepsilon_k}.
\end{equation}
 Define
\[
 \mathbf V_n=
 \begin{cases}
  \text{a fixed enumeration of }W_k,& n=M_k,\\
 \emptyset,& n\notin\{M_k:k\geq1\}.
 \end{cases}
\]
and set $E_{S, A}(\Phi):=\widehat R\cap L^{\prime}(\Phi)$.
Using Lemma~\ref{lem:integrated-blockwise-cost}, \eqref{eq-sec4-4} gives \eqref{eq:integrated-negligible-cost}.
Hence, Proposition~\ref{prop:direct-deterministic-insertions} and Lemma~\ref{Cardi} give
\[\dim_{\rm H} E_{S,A}(\Phi)\ge \frac{1}{\alpha}\]and 
\begin{equation}\label{hougan}E_{S, A}(\Phi)\subset \{x\in E(\Phi) \colon i_n(x)\in S \text{ for all } n\ge 1\}.
\end{equation}
Since $\lVert i_n(x)\rVert_\infty\to \infty$ as $n\to \infty$ for any $x\in E(\Phi)$, by \eqref{hougan} and
Proposition~\ref{upper-d},
\[\dim_{\rm H} E_{S,A}\le \dim_{\rm H}\{x\in E(\Phi) \colon i_n(x)\in S \text{ for all } n\ge 1\}\le\frac{1}{\alpha}.\]
Hence,
\begin{equation}\label{eq-sec4-5}
\dim_{\rm H} E_{S,A}=\dim_{\rm H}\{x\in E(\Phi)\colon \{i_n(x)\colon n\in\mathbb N\}\subset S\}=\frac{1}{\alpha}.
\end{equation}
We now verify the density statement. Put
\[D_{S,A}:=
\bigcup_{n\in\mathbb N}\ \bigcap_{x\in E_{S,A}}\{i_n(x)\}.\]
By construction, every element of each $W_k$ is inserted into the digit sequence of every point of $E_{S,A}(\Phi)$. Hence, for each $k\ge 1,$
\[\bigcup_{j=1}^k W_j \subset D_{S,A}\cap A.\]
Using \eqref{eq-sec4-2}, for each $k\ge 1$ 
\[Q_{N_k}(A\setminus S_{\ast})=\bigcup_{j=1}^k W_j\subset D_{S,A}\cap A\]
and hence
\[\#Q_{N_k}(A\setminus S_{\ast})\le \#Q_{N_k}(D_{S,A}\cap A).\]
Dividing by $\#Q_{N_k}(S)$ and using \eqref{eq-sec4-1}, we obtain
\[\overline{d}(D_{S,A}\cap A | S)
\ge\limsup_{k\to\infty}\frac{\#Q_{N_k}(D_{S,A}\cap A)}{\#Q_{N_k}(S)}\ge
\lim_{k\to\infty}\frac{\#Q_{N_k}(A\setminus S_{\ast})}{\#Q_{N_k}(S)}=
\overline{d}(A| S).\]
Since $D_{S,A}\cap A\subset A$, the reverse inequality
\[\overline{d}(D_{S,A}\cap A | S)\le \overline{d}(A| S)\]
is obvious. Hence,
\begin{equation}\label{eq-sec4-6}
\overline{d}\left(
\bigcup_{n\in\mathbb N} \bigcap_{x\in E_{S,A}}\{i_n(x)\}\cap A \middle| S
\right)=\overline{d}(A| S).
\end{equation}
Combining \eqref{eq-sec4-5} and 
\eqref{eq-sec4-6}, we complete the proof of Theorem~\ref{cor-FS-R}.

\begin{proof}[Proof of Theorem~\ref{cor-FS} (a)]
The proof is obtained by the same construction as that of Theorem~\ref{cor-FS-R}.
Apply Lemma~\ref{find-C} to $I$ and let $I^*\subset I$ be the resulting sparse
subset. Since $I^*$ has zero upper density, we have
$
d(S\setminus I^*)=d(S).$
Choose a sequence $N_k\to\infty$ such that
\[
\frac{\#Q_{N_k}(S\setminus I^*)}{(2N_k+1)^2}
   \longrightarrow d(S),
\]
and define
\[
W_k
=
Q_{N_k}(S\setminus I^*)
\setminus
Q_{N_{k-1}}(S\setminus I^*).
\]
Using $I^*$ for the seed construction and inserting the blocks $W_k$ at
sufficiently sparse levels, exactly as in the proof of Theorem~\ref{cor-FS-R},
Proposition~\ref{prop:direct-deterministic-insertions} and Lemma~\ref{Cardi} give a set $E_S(\Phi)\subset E(\Phi)$ with
\[
\dim_{\rm H} E_S(\Phi)=\dim_{\rm H} E(\Phi)=1.
\]
Since every element of $Q_{N_k}(S\setminus I^*)$ is inserted into the
address of every point of $E_S(\Phi)$, we obtain
\[
d\left(
\bigcup_{n\in\mathbb N}
\bigcap_{x\in E_S(\Phi)}
\{i_n(x)\}\cap S
\right)\ge d(S).
\]
The reverse inequality is immediate, since the set inside the density is
contained in $S$. This proves Theorem~\ref{cor-FS}(a).
\end{proof}


\subsection{Banach-density-preserving insertions and
proof of Theorem~\ref{cor-FS}(b)}
\label{subsec:proof-banach-transference}

We now prove Theorem~\ref{cor-FS}(b).
The strategy is parallel to that used in the proof of
Theorem~\ref{cor-FS-R}.  The main difference is that upper Banach
density is detected on translated squares rather than on centered
squares.  We first choose a sequence of widely separated squares
which asymptotically realizes the upper Banach density.  Their union
will be negligible with respect to centered counting, so that removing
it from the ambient digit set does not affect its polynomial density.

We begin with the following lemma.

\begin{lem}
\label{lem:sparse-banach-blocks}
Let $S\subset\mathbb Z^2.$
Let
\[
 g:[1,\infty)\longrightarrow(0,\infty)
\]
be eventually increasing and satisfy $g(N)\to\infty$ as
$N\to\infty$.  Then there exists a sequence
$(\Pi_k)_{k\geq1}$ of pairwise disjoint axis-parallel integer squares
such that
\begin{equation}
 \text{the side length of \(\Pi_k\) tends to infinity as \(k\to\infty\)},
 \label{eq:banach-side-length}
\end{equation}
\begin{equation}
 \frac{\#(S\cap\Pi_k)}
      {\#(\mathbb Z^2\cap\Pi_k)}
 \longrightarrow d_{\rm B}(S),
 \label{eq:banach-block-density}
\end{equation}
and, with
\[
 U:=
 \mathbb Z^2\cap\bigcup_{k\geq1}\Pi_k,
\]
we have
\begin{equation}
 \#Q_N(U)=o(g(N)).
 \label{eq:sparse-union}
\end{equation}
\end{lem}

\begin{proof}
We first claim the following.
\begin{cla*}
Given $m,R\in\mathbb N$ and
$\eta>0$, there exists an axis-parallel integer square $\Pi$ of side
length $m$ such that
\begin{equation}
 \Pi\cap[-R,R]^2=\emptyset
 \qquad\text{and}\qquad
 \frac{\#(S\cap\Pi)}
      {\#(\mathbb Z^2\cap\Pi)}
 \geq d_{\rm B}(S)-\eta.
 \label{eq:far-dense-square}
\end{equation}
\end{cla*}
\emph{Proof of Claim.}
If $\eta\geq d_{\rm B}(S)$, then $d_{\rm B}(S)-\eta\leq 0$, and any axis-parallel
integer square of side length $m$ disjoint from $[-R,R]^2$ satisfies
\eqref{eq:far-dense-square}. Hence, we may assume that $0<\eta<d_{\rm B}(S)$.

By the definition of the Banach density of $S$, we can choose an axis-parallel
integer rectangle $P$ with arbitrarily large side lengths such that
\[
 \frac{\#(S\cap P)}{\#(\mathbb Z^2\cap P)}
 \geq d_{\rm B}(S)-\frac{\eta}{4}.
\]
Let $\mathcal Q(P,m)$ denote the collection of all axis-parallel integer
squares of side length $m$ contained in $P$, and put
$M:=\#\mathcal Q(P,m).$
Write $P=[a_1,b_1]\times[a_2,b_2]$,
where $a_i,b_i\in\mathbb Z$, and put
$L_i=b_i-a_i$ for $i=1,2$. 
Let $\partial P$ denote the boundary
of $P$ in $\mathbb R^2$, and set
\[
\partial_m P
:=\bigl\{z\in P\cap\mathbb Z^2:
\operatorname{dist}_\infty(z,\partial P)<m\bigr\}.
\]
Here, for $z\in\mathbb R^2$, define
$
\operatorname{dist}_\infty(z,\partial P)
:=\inf_{w\in\partial P}\|z-w\|_\infty.
$
Since $m$ is fixed, we have
$
\#\partial_m P\le C_m(L_1+L_2)
$ and $
M=(L_1-m+1)(L_2-m+1),$ where $C_m$ denotes some constant depending only on $m.$
Consequently,
\begin{align}
\label{sideside}
\frac{\#\partial_m P}{\#(\mathbb Z^2\cap P)}
\longrightarrow 0,
\qquad
\frac{M}{\#(\mathbb Z^2\cap P)}
\longrightarrow 1
\end{align}
as $\min\{L_1,L_2\}\to\infty$.

Every lattice point of $P\setminus\partial_m P$ belongs to exactly
$(m+1)^2$ members of $\mathcal Q(P,m)$. Indeed, for each
$z=(z_1,z_2)\in(P\setminus\partial_m P)\cap\mathbb Z^2$,
the lower-left corner of such a square can be chosen from
$
\{z_1-m,\ldots,z_1\}\times
\{z_2-m,\ldots,z_2\},$
which consists of exactly $(m+1)^2$ lattice points. Therefore, by double counting,
\[
 \sum_{\Pi\in\mathcal Q(P,m)}\#(S\cap\Pi)
 \geq
 (m+1)^2\#\bigl(S\cap(P\setminus\partial_m P)\bigr).
\]
Since every $\Pi\in\mathcal Q(P,m)$ satisfies
$\#(\mathbb Z^2\cap\Pi)=(m+1)^2$, it follows that
\begin{align*}
 \frac1M
 \sum_{\Pi\in\mathcal Q(P,m)}
 \frac{\#(S\cap\Pi)}{\#(\mathbb Z^2\cap\Pi)}
 &\geq
 \frac{\#\bigl(S\cap(P\setminus\partial_m P)\bigr)}{M}\\
 &\geq
 \frac{
 \left(d_{\rm B}(S)-\frac{\eta}{4}\right)
 \#(\mathbb Z^2\cap P)-\#\partial_mP
 }{M}.
\end{align*}
Consequently, by taking the side lengths of $P$ sufficiently large and \eqref{sideside}, we
may arrange that
\begin{equation}
 \frac1M
 \sum_{\Pi\in\mathcal Q(P,m)}
 \frac{\#(S\cap\Pi)}{\#(\mathbb Z^2\cap\Pi)}
 \geq d_{\rm B}(S)-\frac{\eta}{2}.
 \label{eq:far-dense-average}
\end{equation}

On the other hand, the number of axis-parallel integer squares of side
length $m$ which intersect $[-R,R]^2$ is bounded by a constant
$C_{m,R}$ depending only on $m$ and $R$. In particular, this bound is
independent of $P$. Taking $P$ still larger if necessary, we may assume
that
${C_{m,R}}/{M}<{\eta}/{2}.$
Suppose, to the contrary, that every
$\Pi\in\mathcal Q(P,m)$ disjoint from $[-R,R]^2$ satisfies
\[
 \frac{\#(S\cap\Pi)}{\#(\mathbb Z^2\cap\Pi)}
 <d_{\rm B}(S)-\eta.
\]
Since the density of $S$ in any square is at most $1$, and since at most
$C_{m,R}$ squares in $\mathcal Q(P,m)$ intersect $[-R,R]^2$, we obtain
\begin{align*}
 \frac1M
 \sum_{\Pi\in\mathcal Q(P,m)}
 \frac{\#(S\cap\Pi)}{\#(\mathbb Z^2\cap\Pi)}
 &<
 \frac{M-C_{m,R}}{M}(d_{\rm B}(S)-\eta)
 +\frac{C_{m,R}}{M}\\
 &\leq
d_{\rm B}(S)-\eta+\frac{C_{m,R}}{M}\\
 &<
d_{\rm B}(S)-\frac{\eta}{2},
\end{align*}
where we used $\eta<d_{\rm B}(S)\leq 1$ in the second inequality.
This contradicts \eqref{eq:far-dense-average}. Hence there exists
$\Pi\in\mathcal Q(P,m)$ such that
\[
 \Pi\cap[-R,R]^2=\emptyset
 \quad\text{and}\quad
 \frac{\#(S\cap\Pi)}{\#(\mathbb Z^2\cap\Pi)}
 \geq d_{\rm B}(S)-\eta,
\]
which proves \eqref{eq:far-dense-square}.
\qed

Choose integers
$
 1\leq m_1<m_2<\cdots,\qquad m_k\longrightarrow\infty.$
We construct $\Pi_k$ inductively.  Suppose that
$\Pi_1,\ldots,\Pi_{k-1}$ have already been chosen.  Choose $R_k$
sufficiently large so that $[-R_k,R_k]^2$ contains all the preceding
squares and
\begin{equation}
 \sum_{j=1}^{k-1}(m_j+1)^2+(m_k+1)^2
 \leq 2^{-k}g(R_k).
 \label{eq:sparse-placement}
\end{equation}
Using \eqref{eq:far-dense-square}, choose an integer square $\Pi_k$
of side length $m_k$, disjoint from $[-R_k,R_k]^2$, such that
\begin{equation}
 \frac{\#(S\cap\Pi_k)}
      {\#(\mathbb Z^2\cap\Pi_k)}
 \geq d_{\rm B}(S)-\frac1k.
 \label{eq:banach-density-lower}
\end{equation}
The squares are pairwise disjoint and
\eqref{eq:banach-side-length} holds.

Since their side lengths tend to infinity, the definition of the
upper Banach density gives
\[
 \limsup_{k\to\infty}
 \frac{\#(S\cap\Pi_k)}
      {\#(\mathbb Z^2\cap\Pi_k)}
 \leq d_{\rm B}(S).
\]
Together with \eqref{eq:banach-density-lower}, this proves
\eqref{eq:banach-block-density}.

Put
$\rho_k:=
 \min\{\|i\|_\infty:i\in\Pi_k\cap\mathbb Z^2\}.$
Increasing $R_k$ further if necessary during the construction, we may
assume that
$\rho_1<\rho_2<\cdots.$
If
$\rho_k\leq N<\rho_{k+1},$
then only the first $k$ squares can intersect $Q_N$, and hence
\[
 \#Q_N(U)
 \leq
 \sum_{j=1}^k(m_j+1)^2
 \leq
 2^{-k}g(R_k)
 \leq
 2^{-k}g(N).
\]
Since $k\to\infty$ as $N\to\infty$, this proves
\eqref{eq:sparse-union}.
\end{proof}


\begin{proof}[Proof of Theorem~\ref{cor-FS}(b)]
Let $S\subset I$ satisfy
$
 d_{\rm B}(S)>0.
$
By assumption, $I$ has polynomial density with exponent $1$.
Thus there exist $\beta\geq0$ and $\gamma_1,\gamma_2>0$ such that,
for all sufficiently large $N$,
\begin{equation}
 \gamma_1\frac{N^2}{(\log N)^\beta}
 \leq
 \#Q_N(I)
 \leq
 \gamma_2\frac{N^2}{(\log N)^\beta}.
 \label{eq:I-polynomial-density}
\end{equation}
Apply Lemma~\ref{lem:sparse-banach-blocks} with
\[
 g(N):=\frac{N^2}{(\log N)^\beta}.
\]
Let $(\Pi_k)_{k\geq1}$ be the resulting sequence of squares and put
\begin{equation}
 W_k:=S\cap\Pi_k,
 \qquad
 W:=\bigcup_{k\geq1}W_k.
 \label{eq:Wk-banach}
\end{equation}
The sets $W_k$ are finite and pairwise disjoint, and
\begin{equation}
 \frac{\#W_k}{\#(\mathbb Z^2\cap\Pi_k)}
 \longrightarrow d_{\rm B}(S).
 \label{eq:Wk-density}
\end{equation}
Moreover,
\eqref{eq:sparse-union} and
\eqref{eq:I-polynomial-density} imply
\begin{equation}
 \frac{\#Q_N(W)}{\#Q_N(I)}
 \longrightarrow0.
 \label{eq:W-zero-relative-density}
\end{equation}
Set
$B_0:=I\setminus W.$
It follows from \eqref{eq:W-zero-relative-density} that
\[
 \#Q_N(B_0)
 =
 (1+o(1))\#Q_N(I).
\]
In particular, $B_0$ also has polynomial density with exponent $1$.

We use $B_0$ directly for the seed construction.  Since $B_0$
has polynomial density with exponent $1$, there exist
$\widetilde\gamma_1,\widetilde\gamma_2>0$ such that, for all
sufficiently large $N$,
\[
 \widetilde\gamma_1\frac{N^2}{(\log N)^\beta}
 \leq \#Q_N(B_0)
 \leq
 \widetilde\gamma_2\frac{N^2}{(\log N)^\beta}.
\]
Choose integers $L>2$ and $t>L$ sufficiently large that
\[
 c_0:=2^{-\beta}\widetilde\gamma_1L^2
      -\widetilde\gamma_2>0,
\]
and define
\[
 I^{(n)}
 =Q_{Lt^n}(B_0)\setminus Q_{t^n-1}(B_0),
 \qquad n\geq1.
\]
Since $t>L$, these alphabets are pairwise disjoint, and
\[
 t^n\leq\|i\|_\infty\leq Lt^n
 \qquad(i\in I^{(n)}).
\]
Moreover, for all $n\geq1$,
\begin{align}
 \#I^{(n)}
 &\geq \#Q_{Lt^n}(B_0)-\#Q_{t^n}(B_0)
 \notag\\
 &\geq
 \frac{\widetilde\gamma_1L^2t^{2n}}
      {(n\log t+\log L)^\beta}
 -\frac{\widetilde\gamma_2t^{2n}}
       {(n\log t)^\beta}
 \notag\\
 &\geq \frac{c_0t^{2n}}{(n\log t)^\beta},
 \label{eq:banach-seed-cardinality}
\end{align}
where the last inequality uses $L<t$.
By increasing $t$ if necessary, we may assume that
$I^{(n)}$ is non-empty for every $n\geq1$.

Let $D_*$ be as in \eqref{eq:integrated-seed-dimension}.
By \eqref{eq:banach-seed-cardinality},
$
 \sum_{j=1}^n\log\#I^{(j)}
 \geq n(n+1)\log t-O(n\log n),$
and hence $D_*\geq1$.  On the other hand,
$\#I^{(n)}\leq(2Lt^n+1)^2$ implies $D_*\leq1$.
Consequently,
\begin{equation}
 D_*=1.
 \label{eq:banach-seed-D}
\end{equation}

Since $B_0\cap W=\emptyset$, we can insert all the
elements of $W_k$ into the seed construction at sufficiently sparse
levels.

Let $(\varepsilon_k)_{k\geq1}$ be a decreasing sequence of positive
numbers converging to $0$.  We inductively choose
$0=M_0<M_1<M_2<\cdots$
as follows.  Suppose that $M_0,\ldots,M_{k-1}$ have already been
chosen.  Since
\[
 \prod_{j=M_{k-1}+1}^{M}
 \frac1{t^{2j}}
 \longrightarrow0
 \qquad(M\to\infty),
\]
we may choose $M_k>M_{k-1}$ sufficiently large that
\begin{equation}
 \prod_{i\in W_k}\frac1{\|i\|_\infty^2}
 \geq
 \left(
 \prod_{j=M_{k-1}+1}^{M_k}\frac1{t^{2j}}
 \right)^{\varepsilon_k}.
 \label{eq:banach-insertion-product}
\end{equation}
Define
\[
 V_n=
 \begin{cases}
 \text{a fixed enumeration of }W_k,
     & n=M_k,\\[2mm]
 \emptyset,
     & n\notin\{M_k:k\geq1\}.
 \end{cases}
\]
Since the sets $W_k$ are pairwise disjoint and
\[
 B_0\cap W=\emptyset,
\]
all inserted symbols are pairwise distinct and none of them belongs
to any seed alphabet $I^{(n)}$.

Let $\widehat R$ be the inserted Moran set associated with
$(I^{(n)})_{n\geq1}$ and $(V_n)_{n\geq1}$ as in
Section~3, and put
\begin{equation}
 F_S(\Phi):=\widehat R\cap L'(\Phi).
 \label{eq:FS-definition}
\end{equation}

By
\eqref{eq:banach-insertion-product} and
Lemma~\ref{lem:integrated-blockwise-cost}, the cumulative insertion
cost satisfies
$\Lambda_n=o(n^2).$
Therefore Proposition~\ref{prop:direct-deterministic-insertions} and
\eqref{eq:banach-seed-D} give
\begin{equation}
 \dimH F_S(\Phi)\geq1,
 \qquad
 F_S(\Phi)\subset E(\Phi).
 \label{eq:FS-lower-dimension}
\end{equation}

On the other hand, for every $x\in E(\Phi)$ the digits
$i_n(x)$ are pairwise distinct.  Since every bounded subset of
$\mathbb Z^2$ is finite, this implies
$\|i_n(x)\|_\infty\longrightarrow\infty.$
Hence,
$
 E(\Phi)\subset G_\Phi(I).$
Proposition~\ref{upper-d} and the polynomial density of $I$ with exponent $1$
imply
$\dimH E(\Phi)\leq1.$
Combining this with
\eqref{eq:FS-lower-dimension}, we obtain
\begin{equation}
 \dimH F_S(\Phi)
 =
 \dimH E(\Phi)
 =
 1.
 \label{eq:FS-full-dimension}
\end{equation}

It remains to verify the upper Banach density statement.
Define the set of digits common to all points of $F_S(\Phi)$ by
\begin{equation}
 D(F_S)
 :=
 \bigcup_{n\geq1}
 \bigcap_{x\in F_S(\Phi)}
 \{i_n(x)\}.
 \label{eq:common-digit-set}
\end{equation}
Every element of $W_k$ is inserted at a prescribed position in the
address of every point of $F_S(\Phi)$.  Since points of $F_S(\Phi)$
have unique addresses,
$W_k\subset D(F_S)
 $ for each $k\ge 1.$

Moreover, every digit appearing in the construction belongs either
to $B_0$ or to $W$.  Since
$B_0\cap W=\emptyset,
 W_k=S\cap\Pi_k,$
and the squares $\Pi_k$ are pairwise disjoint, we have
\begin{equation}
 D(F_S)\cap S\cap\Pi_k=W_k
 \qquad(k\geq1).
 \label{eq:common-digits-on-block}
\end{equation}
Using \eqref{eq:Wk-density}, we obtain
\begin{align*}
 d_{\rm B}(D(F_S)\cap S)
 &\geq
 \lim_{k\to\infty}
 \frac{\#(D(F_S)\cap S\cap\Pi_k)}
      {\#(\mathbb Z^2\cap\Pi_k)}
 \\
 &=
 \lim_{k\to\infty}
 \frac{\#W_k}{\#(\mathbb Z^2\cap\Pi_k)}=d_{\rm B}(S).
\end{align*}
Since
$D(F_S)\cap S\subset S,$
the reverse inequality is immediate.  Thus
\begin{equation}
 d_{\rm B}
 \left(
 \bigcup_{n\geq1}
 \bigcap_{x\in F_S(\Phi)}
 \{i_n(x)\}
 \cap S
 \right)
 =
 d_{\rm B}(S).
 \label{eq:FS-banach-density}
\end{equation}
Together with \eqref{eq:FS-full-dimension}, this proves
Theorem~\ref{cor-FS}(b).
\end{proof}

\section{Applications}
Since the work of Hurwitz~\cite{H87}, many continued-fraction
algorithms for complex numbers have been introduced; see, for example,
\cite{DN,H1902,KST76,Lak73,LeV52,SKT75}.
In this final section, we consider three representative complex
continued-fraction systems to illustrate applications of our main
results.

\if0\subsection{Background}

What the authors call `upper density' in \cite{FK78} appears to be an upper Banach density, see \cite[p.276]{FK78}.
In \cite{BL},  `having positive upper Banach density' is defined on p.725
but `upper Banach density' itself is not defined.

In (\url{https://mathlog.info/articles/PxvytNBwCOsonA3XtIx6}) they defined the upper Banach density of $S\subset\mathbb Z^k$ as follows:

\[d_{\rm B}'(S)=\limsup_{N\to\infty}\frac{1}{(2N+1)^k}\sup_{M\in\mathbb Z^k}\#(S\cap (M+[-N,N]^k)).\]
\textcolor{red}{Which of the definitions of the Banach density are you using?}\textcolor{cyan}{I've not treated Banach density yet.}
I am not sure if $d_{\rm B}(S)$ coincides with $d_{\rm B}'(S)$. It is easy to see $d_{\rm B}(S)\geq d_{\rm B}'(S)$. If $k=1$ then $d_{\rm B}(S)= d_{\rm B}'(S)$ for any $S\subset\mathbb Z$. I suspect the equality does not hold in general if $k\geq2$.

\fi

\subsection{The Hurwitz continued fraction}Hurwitz~\cite{H87} introduced a nearest-integer-type continued
fraction expansion for complex numbers based on Gaussian integers. Various metric and fractal properties of the Hurwitz continued fraction
have subsequently been studied; see, for example,
\cite{BGH25, Ger,NT2} and the references therein.
We recall this expansion below. 
Under the natural identification of 
 $\mathbb C$ with $\mathbb R^2$, the Gaussian field $\mathbb Q(\sqrt{-1})$ and the ring $\mathbb Z[\sqrt{-1}]$ of Gaussian integers are identified  with $\mathbb Q^2$ and $\mathbb Z^2$ respectively.
Define \[U=\left\{x\in \mathbb C\colon  
-\frac{1}{2}\le {\rm Re}(x)<\frac{1}{2},\ -\frac{1}{2}\leq {\rm Im}(x)< \frac{1}{2}\right\}.\]
For $x\in \mathbb C$ we define \[\lfloor x\rfloor=\left\lfloor {\rm Re}(x)+\frac{1}{2}\right\rfloor+\sqrt{-1} \left\lfloor {\rm Im}(x)+\frac{1}{2}\right\rfloor,\]
which is one of the at most four Gaussian integers that are nearest to $x$ in the Euclidean metric.
Define $H\colon U\setminus\{0\}\rightarrow U$ by \begin{equation}\label{H-map}H(x)=\frac{1}{x}-\left\lfloor \frac{1}{x}\right\rfloor.\end{equation}
For each $x\in U\setminus\{0\}$, define a sequence $i_1(x)$, $i_2(x),\ldots$ of nonzero Gaussian integers 
inductively by 
\[i_{n}(x)=\left\lfloor \frac{1}{H^{n-1}(x)}\right\rfloor\text{ for }n=1,2,\ldots.\]
If $H^{n}(x)=0$ for some $n\geq1$ then $i_{n+1}(x)$ is not defined.
We have $H^n(x)\neq0$ for every $n\geq1$
if and only if $x\in U\setminus\mathbb Q(\sqrt{-1})$, and in this case $x$
   has the unique infinite expansion of the form 
\[
x=
\cfrac{1}{
k_1(x)+l_1(x)\sqrt{-1}
+\cfrac{1}{
k_2(x)+l_2(x)\sqrt{-1}
+\cfrac{1}{
k_3(x)+l_3(x)\sqrt{-1}
+\ddots
}
}
}.
\]
where $i_n(x)=k_{n}(x)+l_n(x)\sqrt{-1}$.
See \cite[Theorem~6.1]{DN} for example.
This expansion is called the
 {\it Hurwitz continued fraction expansion} \cite{H87}.

\begin{figure}
\begin{center}
\centering
\includegraphics[width=9cm,trim=0 2cm 0 0,clip]{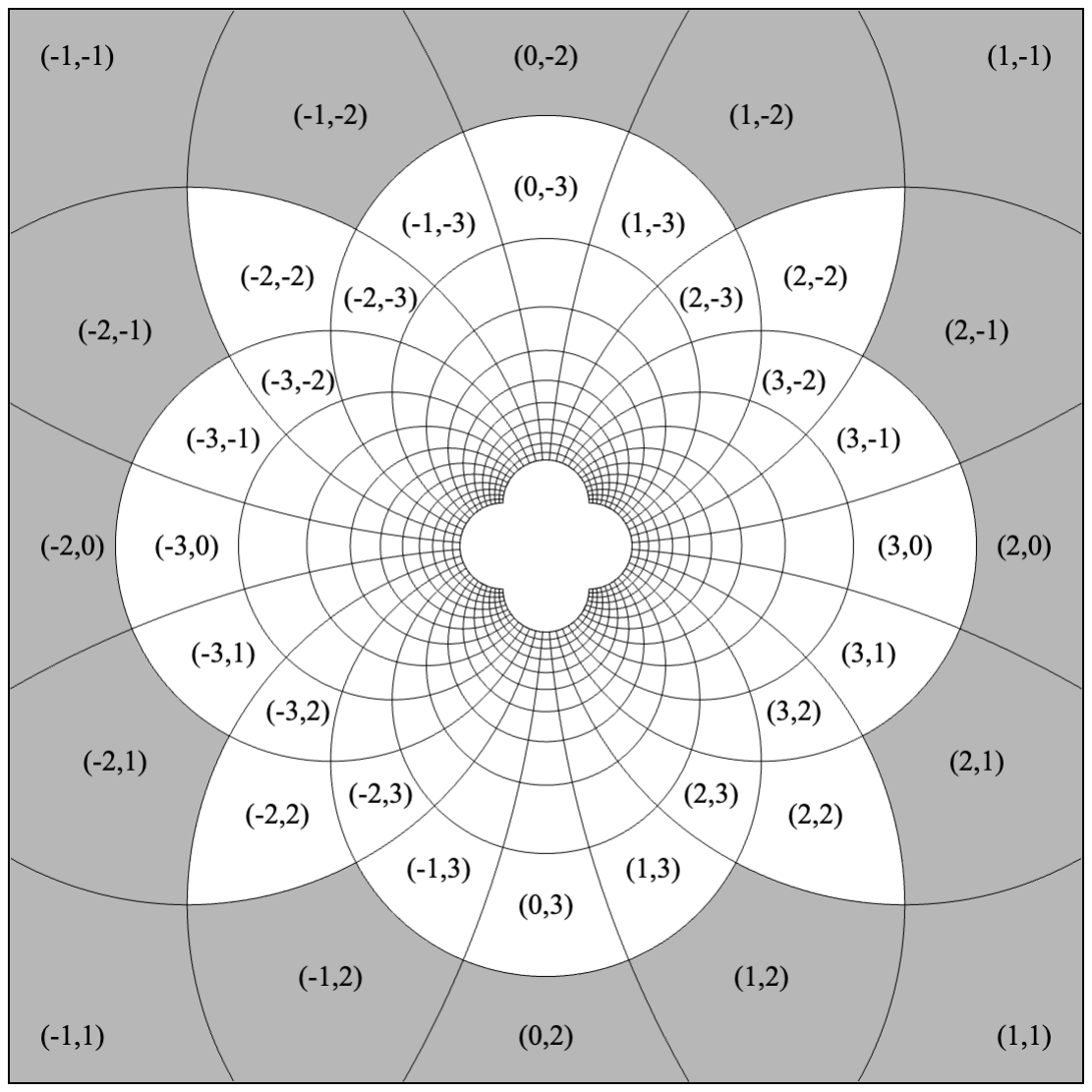}
\vspace{-30pt}
\caption
{The set $U\setminus\{0\}$ is partitioned into $U_i\cap U$,  $i\in\mathbb D_1$.
The shaded region indicates $\bigcup\{U_{i}\cap U\colon i\in\mathbb D_1\setminus\mathbb D_2\}$. The limit set of the infinite conformal IFS $\{\phi_i\colon i\in\mathbb D_2\}$ on $\Delta=\overline{U}$ associated with the Hurwitz continued fraction is contained in the unshaded region.}\label{fig3}
\end{center}
\end{figure}

There is a canonical infinite iterated function system 
 that represents the regular continued fraction,  whereas this is not the case for the Hurwitz continued fraction.
However, neglecting finitely many cylinders one can extract an infinite iterated function system from the Hurwitz continued fraction. 
Put 
\[\mathbb D_{1}=\{(k,l)\in\mathbb Z^2\colon k^2+l^2\geq 2\}\ \text{  
and }\
\mathbb D_{2}=\{(k,l)\in\mathbb Z^2\colon k^2+l^2\geq 8\}.\]
For each 
$i\in\mathbb D_1$, 
define a set
\[U_i=\left\{x\in\mathbb C\colon\frac{1}{x}\in i+U\right\}.\]
It is easy to see that: 
\begin{itemize}  
\item 
$U_{i}\cap \partial U\neq \emptyset$
if $i\in\mathbb D_1\setminus\mathbb D_2$ and
$\overline{U_{i}}\subset U$ if $i\in\mathbb D_2$;
\item  
$\bigcup\{U_{i}\cap U\colon i\in\mathbb D_1\}=U\setminus\{0\}$.
\end{itemize}

Put $\Delta=\overline{U}$, and 
 for each $i=(k,l)\in\mathbb D_2$ 
 define a map $\phi_{i}\colon  \Delta\to\mathbb C$ by
\[\phi_{i}(x)= \frac{1}{x+k+l\sqrt{-1}}.\]
Note that $\phi_{i}|_{{\rm int}U}$ is univalent, $\phi_{i}(U)=U_{i}$ and $\phi_{i}(\Delta)\subset \Delta$.
The estimates in \cite[Lemma~4.2]{NT2} give (A1)--(A3)
and \eqref{CF-der}. Since  $|z|\le1/\sqrt2$ for $z\in\Delta$ and
$|i|\ge\sqrt8$ for $i\in\mathbb D_2$, we have
$
|\phi_i(z)|\le{\sqrt2}/{3}<1/2.
$
Hence, $\phi_i(\Delta)\subset\operatorname{int}\Delta$ for every
$i\in\mathbb D_2$, and (A1) gives
$\mathcal O(\Phi)=\emptyset$, where
$\Phi=\{\phi_i:i\in\mathbb D_2\}$.
Thus (A4) holds, and $\Phi$ is a $2$-decaying conformal IFS
on $\Delta$.
See \textsc{Figure}~\ref{fig3}.

Consider the set
\[E_{\rm H}=\{x\in U\setminus\mathbb Q^2\colon i_m(x)\neq i_{n}(x)\ \text{ for all } m,n\in\mathbb N\text{ with }m\neq n\}.\]
Since $\mathbb Z^2\setminus D_2$ is finite, for every
$S\subset\mathbb Z^2$ we have
$d_{\rm B}(S\cap D_2)=d_{\rm B}(S).$
Moreover, $\mathbb P^2\subset D_2$.
Hence Theorems~\ref{cor-N} and~\ref{GT-trans}, applied to
$\{\phi_i:i\in D_2\}$, yield the following statements.
\begin{thm}\ 
\begin{itemize}
\item[(a)] For any 
$S\subset\mathbb Z^2$ with $d_{\rm B}(S)>0$, 
there exists $F_S\subset E_{\rm H}$ with Hausdorff dimension 
$1$ such that
  for any $h_1,\ldots,h_\ell\in\mathbb Q_0[X]$, $\ell\geq2$, $v_1,\ldots,v_\ell\in \mathbb Z^2$  
  there are $n\in\mathbb N$, $u\in\mathbb Z^2$ and 
  $n_1,\ldots,n_\ell\in \mathbb N$
  such that  \[\{u+h_j(n)v_j\colon 1\leq j\leq \ell\}\subset S\]
  and 
  for any $x\in F_S$, \[i_{n_j}(x)=u+h_j(n)v_j\ \text{ for }1\leq j\leq\ell.\]

\item[(b)] For any $A\subset\mathbb P^2$ with $\overline{d}(A|\mathbb P^2)>0$, there exists $F_{\mathbb P^2,A}\subset E_{\rm H}$ with Hausdorff dimension $1$ 
  such that for every $x\in F_{\mathbb P^2,A}$,
\[
\{i_n(x):n\in\mathbb N\}\subset\mathbb P^2,
\]
and
\[
\{i_n(x):n\in\mathbb N\}\cap A
\]
contains infinitely many homothetic copies of any finite subset of
$\mathbb Z^2$.
  \end{itemize}
\end{thm}

\begin{cor}\label{cor-Hur}
The set
\[\{x\in E_{\rm H}\colon i_n(x)\in\mathbb P^2\ \text{ for all }n\geq1\}\] 
is of Hausdorff dimension $1$.\end{cor}

The Hurwitz system above was previously used in~\cite{NT3} to show
that the set of points whose partial quotients contain infinitely many
homothetic copies of every finite subset of $\mathbb Z^2$ has
Hausdorff dimension $1$. The results above strengthen and place this
phenomenon in the general framework of the present paper by allowing
the patterns to be prescribed inside arbitrary positive-density digit
sets.

\subsection{Kaneiwa--Shiokawa--Tamura's continued fraction}
\label{KST-sec}

Let
\[
    \zeta=\frac{1+\sqrt{-3}}{2},
    \qquad
    \overline{\zeta}=\frac{1-\sqrt{-3}}{2},
\]
and put
\[
    V=\{s\zeta+t\overline{\zeta}\in\mathbb C:
    0\leq s<1,\ 0\leq t<1\}.
\]
The ring of integers of $\mathbb Q(\sqrt{-3})$ is
$\mathcal O_{-3}
    =
    \{m\zeta+n\overline{\zeta}:m,n\in\mathbb Z\}.$
Put
$\Delta=\overline{V}.$

We recall the complex continued fraction algorithm introduced by
Shiokawa, Kaneiwa and Tamura~\cite{SKT75}.
Every $x\in\mathbb C$ can be written uniquely in the form
$x=s\zeta+t\overline{\zeta},
    s,t\in\mathbb R.$
Define
\[
    \lfloor x\rfloor_{\rm KST}
    =
    \lfloor s\rfloor\zeta+\lfloor t\rfloor\overline{\zeta}.
\]
Equivalently, $\lfloor x\rfloor_{\rm KST}$ is the unique element of
$\mathcal O_{-3}$ such that
$x\in \lfloor x\rfloor_{\rm KST}+V.$
The KST map is then defined by
\[
    T_{\rm KST}:V\setminus\{0\}\longrightarrow V,
    \qquad
    T_{\rm KST}(x)
    =
    \frac1x-\left\lfloor\frac1x\right\rfloor_{\rm KST}.
    \tag{5.2}
\]
For $x\in V\setminus\{0\}$, define $k_n(x),l_n(x)\in\mathbb Z$ by
\[
    \left\lfloor
    \frac{1}{T_{\rm KST}^{\,n-1}(x)}
    \right\rfloor_{\rm KST}
    =
    k_n(x)\zeta+l_n(x)\overline{\zeta},
    \qquad n\geq1,
\]
whenever $T_{\rm KST}^{\,n-1}(x)\neq0$, and put
$ i_n(x):=(k_n(x),l_n(x)).$
As shown in~\cite{SKT75}, the algorithm terminates after finitely many
steps if and only if $x\in\mathbb Q(\sqrt{-3})$. Hence, if
$x\in V\setminus\mathbb Q(\sqrt{-3})$, the sequence
$(i_n(x))_{n\geq1}$ is defined for all $n$, and
\[
    x=
    \cfrac{1}{
      k_1(x)\zeta+l_1(x)\overline{\zeta}
      +\cfrac{1}{
        k_2(x)\zeta+l_2(x)\overline{\zeta}
        +\cfrac{1}{
          k_3(x)\zeta+l_3(x)\overline{\zeta}+\ddots
        }
      }
    }.
\]
\begin{figure}
\begin{center}
\includegraphics[]{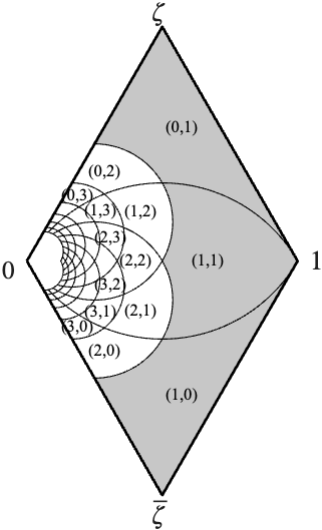}
\caption
{The set $V\setminus\{0\}$ is partitioned into $V_i\cap V$,  $i\in\mathbb N_0^2\setminus \{(0, 0)\}$.
The shaded region indicates $\bigcup\{V_{i}\cap V\colon i\in\mathbb N_0^2\setminus I_{\rm KST}\}$. The limit set of the infinite conformal IFS $\{\phi_i\colon i\in I_{\rm KST}\}$ on $\Delta=\overline{V}$ associated with the KST continued fraction is contained in the unshaded region.}\label{figKST}
\end{center}
\end{figure}
We next describe an infinite conformal IFS naturally associated with
this expansion. 
Put $\mathbb N_0:=\mathbb N\cup \{0\}.$ For
$i=(i_1,i_2)\in\mathbb N_0^2\setminus\{(0,0)\},$
define
\[
    V_i
    =
    \left\{
       x\in\mathbb C:
       \frac1x\in i_1\zeta+i_2\overline{\zeta}+V
    \right\}.
\]
A direct computation gives  
\[
    V_i\subset V
    \quad\Longleftrightarrow\quad
    i\in I_{\rm KST}:=\mathbb N_0^2\setminus\{(0, 0), (0,1), (1,0), (1, 1)\}.
\]
For each $i=(k,l)\in I_{\rm KST}$, define
$\phi_i:\Delta \to \Delta$
by
\[
    \phi_i(x)
    =
    \frac{1}{x+k\zeta+l\overline{\zeta}}.
\]
Then
$
\phi_i(V)=V_i,
\phi_i(\Delta)=\overline{V_i}.$
Set
$
    \Phi_{\rm KST}
    =
    \{\phi_i:i\in I_{\rm KST}\}.$
As in Section~5.1, $\Phi_{\rm KST}$ satisfies
(A1)--(A3) and \eqref{CF-der}. We now verify (A4). For $z=s\zeta+t\overline\zeta\in\Delta$ and $(k,l)\in I_{\rm KST}$,
write $u=k+s$, $v=l+t$ and $D=u^2-uv+v^2$.
Then
\[
\phi_{k,l}(z)
=\frac{v}{D}\zeta+\frac{u}{D}\overline\zeta.
\]
For $(k,l)\in I_{\rm KST}$,
$M:=\max\{u,v\}\ge2$ and
$D\ge 3M^2/4>M$. Then
$\phi_{k,l}(\Delta)\subset V$, and hence
$L(\Phi_{\rm KST})\subset V$. For $i=(k,l)\in I_{\rm KST}$, put
$
a_i:=k\zeta+l\bar\zeta.$
If $\phi_i(y)=\phi_j(y')$ for $y,y'\in L(\Phi_{\rm KST})$,
then $a_i+y=a_j+y'$. Since the translates $a_i+V$ are
pairwise disjoint, we must have $i=j$. Thus
$\mathcal O(\Phi_{\rm KST})=\emptyset$, and (A4) holds.
Hence, $\Phi_{\rm KST}$ is a $2$-decaying conformal IFS.
See \textsc{Figure}~\ref{figKST}.

%

Consider the set
\[
    E_{\rm KST}
    =
    \left\{
       x\in V\setminus\mathbb Q(\sqrt{-3}):
       i_m(x)\neq i_n(x)
       \text{ for all }m,n\in\mathbb N,\ m\neq n
    \right\}.
\]
Notice that, for every $S\subset\mathbb N^2$, removing the single
point $(1,1)$ does not affect upper Banach density:
$d_{\rm B}(S\cap I_{\rm KST})=d_{\rm B}(S).$
Moreover,
$ \mathbb P^2\subset I_{\rm KST}.$
Theorems~\ref{cor-N} and~\ref{GT-trans}, applied to $\Phi_{\rm KST}$, yield the
following statements.

\begin{thm}
\label{thm-KST}\
\begin{enumerate}
\item[(a)]
For any $S\subset\mathbb N^2$ with $d_{\rm B}(S)>0$, there exists
$F_S\subset E_{\rm KST}$ with Hausdorff dimension $1$ such that,
for any
$h_1,\ldots,h_\ell\in\mathbb Q_0[X],    v_1,\ldots,v_\ell\in\mathbb Z^2,
    \ell\geq2,$
there exist
$n\in\mathbb N,
    u\in\mathbb  Z^2,
    n_1,\ldots,n_\ell\in\mathbb N$
such that
\[
    \{u+h_j(n)v_j:1\leq j\leq\ell\}\subset S,
\]
and, for every $x\in F_S$,
\[
    i_{n_j}(x)=u+h_j(n)v_j \ \text{for}\ 1\leq j\leq\ell.
\]

\item[(b)]
For any $A\subset\mathbb P^2$ with
$d(A\mid\mathbb P^2)>0,$
there exists
$F_{\mathbb P^2,A}\subset E_{\rm KST}$
with Hausdorff dimension $1$ such that, for every
$x\in F_{\mathbb P^2,A}$,
\[
    \{i_n(x):n\in\mathbb N\}\subset\mathbb P^2,
\]
and
\[
    \{i_n(x):n\in\mathbb N\}\cap A
\]
contains infinitely many homothetic copies of any finite subset of
$\mathbb Z^2$.
\end{enumerate}
\end{thm}

As an immediate consequence of Theorem~\ref{thm-KST}(b), we obtain
the following.

\begin{cor}
\label{cor-KST}
The set
\[
    \left\{
       x\in E_{\rm KST}:
       i_n(x)\in\mathbb P^2
       \text{ for every }n\geq1
    \right\}
\]
has Hausdorff dimension $1$.
\end{cor}

\subsection{Mauldin-Urba\'nski's infinite conformal IFS}
\begin{figure}
\begin{center}
\includegraphics[]{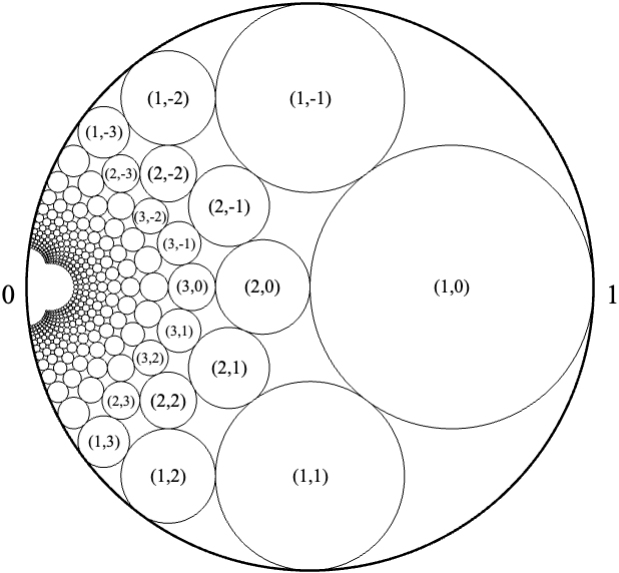}
\caption
{Mauldin--Urba\'nski's infinite conformal IFS
$\Phi_{\rm MU}=\{\phi_i:i\in\mathbb N\times\mathbb Z\}$ on
$\Delta=\{z\in\mathbb C:|z-1/2|\leq1/2\}.$
Although $L(\Phi_{\rm MU})\subset\operatorname{int}\Delta$,
its closure meets $\partial\Delta$; in particular, the limit set accumulates at $0$.}
\label{figMU}
\end{center}
\end{figure}

The final example has played an important role in the creation of dimension theory of conformal IFSs. 
Let 
\[\Delta=\left\{x\in\mathbb C\colon\left|x-\frac{1}{2}\right|\leq\frac{1}{2}\right\}.\]
For each $i=(k,l)\in\mathbb N\times\mathbb Z$, define a map $\phi_{i}\colon\Delta\to\Delta$ by 
\[\phi_{i}(x)=\frac{1}{x+k+l\sqrt{-1}}.\]
The infinite conformal IFS $\Phi_{\rm MU}=\{\phi_{i}\colon i\in\mathbb N\times\mathbb Z\}$ was introduced by 
Mauldin and Urba\'nski \cite[Section~6]{MU96}. Generalizations of this system to families of conformal IFSs associated
with complex continued fractions have been studied by Inui, Okada and
Sumi~\cite{IOS}, and by Inui and Sumi~\cite{IS}, from the viewpoints
of Hausdorff dimension and packing measure, respectively. As in Section~5.1, $\Phi_{\rm MU}$ satisfies
(A1)--(A3) and \eqref{CF-der}. 
We verify that (A4) holds.
Since every $x\in L(\Phi_{\rm MU})$ belongs to $\Delta\setminus\{0\}$, $\operatorname{Re}x>0$.
For $y=\phi_{k,l}(x)$ we have
$\operatorname{Re}(1/y)=k+\operatorname{Re}x>1,$
which is equivalent to $y\in\operatorname{int}\Delta$.
Thus $L(\Phi_{\rm MU})\subset\operatorname{int}\Delta$. By (A1),
$\mathcal O(\Phi_{\rm MU})=\emptyset$. Hence, $\Phi_{\rm MU}$ is a $2$-decaying conformal IFS.
See \textsc{Figure}~\ref{figMU}.
For any $x\in L(\Phi_{\rm MU})$ we have 
\[
x=
\cfrac{1}{
k_1(x)+l_1(x)\sqrt{-1}
+\cfrac{1}{
k_2(x)+l_2(x)\sqrt{-1}
+\cfrac{1}{
k_3(x)+l_3(x)\sqrt{-1}
+\ddots
}
}
}.
\]
Here, $
i_n(x)=(k_n(x),l_n(x)).$

Mauldin and Urba\'nski showed that $L(\Phi_{\rm MU})$ has an interesting geometric property not present in finite conformal IFSs: 
writing $h=\dim_{\rm H} L(\Phi_{\rm MU})$, 
the $h$-dimensional Hausdorff measure of $L(\Phi_{\rm MU})$ is zero,
whereas its $h$-dimensional packing measure is positive and finite.


Consider the set
\[E_{\rm MU}=\{x\in L^{\prime}(\Phi_{\rm MU})\colon i_m(x)\neq i_{n}(x)\ \text{ for all } m,n\in\mathbb N\text{ with }m\neq n\}.\]
The index set $\mathbb N\times\mathbb Z$ has polynomial density
with exponent $1$, and
$\mathbb P^2\subset\mathbb N\times\mathbb Z.$
Hence, applying Theorems~\ref{cor-N} and~\ref{GT-trans} to
$\Phi_{\rm MU}$ yields the following.
\begin{thm}\
\begin{itemize}
\item[(a)] For any $S\subset \mathbb N\times\mathbb Z$ with $d_{\rm B}(S)>0$, there exists $F_S\subset E_{\rm MU}$ with Hausdorff dimension $1$  such that
  for any $h_1,\ldots,h_\ell\in\mathbb Q_0[X]$, $\ell\geq2$, $v_1,\ldots,v_\ell\in \mathbb Z^2$  
  there are $n\in\mathbb N$, $u\in\mathbb Z^2$ and
  $n_1,\ldots,n_\ell\in \mathbb N$
  such that  \[\{u+h_j(n)v_j\colon 1\leq j\leq \ell\}\subset S,\]
  and 
  for any $x\in F_S$, \[i_{n_j}(x)=u+h_j(n)v_j\ \text{ for }1\leq j\leq\ell.\]

\item[(b)] For any $A\subset\mathbb P^2$ with $\overline{d}(A|\mathbb P^2)>0$, there exists $F_{\mathbb P^2,A}\subset E_{\rm MU}$ with Hausdorff dimension $1$ such that
  for every $x\in F_{\mathbb P^2,A}$,
\[
\{i_n(x):n\in\mathbb N\}\subset\mathbb P^2,
\]
and
\[
\{i_n(x):n\in\mathbb N\}\cap A
\]
contains infinitely many homothetic copies of any finite subset of
$\mathbb Z^2$.
  \end{itemize}
\end{thm}

\subsection*{Acknowledgments} 
YN was supported by the JSPS KAKENHI 25K17282, Grant-in-Aid for Early-Career Scientists. HT was supported by the
JSPS KAKENHI 25K21999, Grant-in-Aid for Challenging Research (Exploratory) and 26H02003, Grant-in-Aid for Scientific Research (A).
\subsection*{AI Declaration}
During the preparation of this manuscript, the authors used ChatGPT
(OpenAI, GPT-5.6 Sol) to improve the clarity and readability of the English text, to assist in generating MATLAB code for the figures in this paper, 
to assist with mathematical brainstorming, and to refine the
presentation of several arguments. In particular, it was used to explore a refinement of the dimension estimate in Section 3 that allowed us to remove a technical assumption, and to develop an approach to the proof of the assertion concerning upper Banach
density in Theorem~\ref{cor-FS}(b). All mathematical arguments suggested in this process were independently checked, modified where necessary, and incorporated into the manuscript by the authors. The authors take full responsibility for the correctness and content of the paper.

\end{document}